\documentclass[12pt]{amsart}
\usepackage{amssymb,amscd}

\newtheorem{thm}{Theorem}[section]
\newtheorem{prop}[thm]{Proposition}
\newtheorem{lem}[thm]{Lemma}
\newtheorem{rem}[thm]{Remark}
\newtheorem{cor}[thm]{Corollary}
\numberwithin{equation}{section}
\def\Ker{\operatorname{Ker}}
\def\Im{\operatorname{Im}}
\def\Hom{\operatorname{Hom}} 
\newcommand{\Q}{\mathbb{Q}}
\newcommand{\Z}{\mathbb{Z}}
\newcommand{\C}{\mathbb{C}}

\begin{document}

\title{ The integral cohomology ring of $E_{8}/T^1 \! \cdot \! E_7$}
\author{Masaki Nakagawa$^*$ }
\pagestyle{plain}



\thanks{
$^*$ Partially supported by the Grant-in-Aid for Scientific  Research 
(C)  21540104,   Japan Society for the Promotion of Science. 
}

\address{Department of General Education \endgraf
                      Kagawa National College of Technology \endgraf
                      Takamatsu 761-8058 \\ Japan}
\email{nakagawa@t.kagawa-nct.ac.jp}

\maketitle 
 
\begin{abstract}
    We determine the integral cohomology ring of the homogeneous space
  $E_8/T^1 \! \cdot \! E_7$ by the Borel presentation  and a method due to 
  Toda. Then using the Gysin exact sequence associated with the circle bundle $S^1 \rightarrow 
  E_8/ E_7 \rightarrow E_8/T^1 \! \cdot \! E_7$, we also determine the integral cohomology
  of $E_8/E_7$. 
\end{abstract} 

\section{Introduction}
Let $G$ be a compact connected Lie group and $H$ a centralizer of a toral subgroup. Then the homogeneous space
$G/H$ is called a {\it generalized flag manifold}.\footnote{Generalized flag manifolds exhaust all compact, simply connected 
homogeneous K\"{a}hler manifolds as was shown by  Wang \cite{Wang54} (see also \cite{Bor54}, \cite[13.5]{Bor-Hir58}). 
For this reason, they are 
sometimes called {\it K\"{a}hlerian C-spaces}.}  It is a natural generalization of a (full or complete) flag manifold $G/T$, where 
$T$ is a maximal torus of $G$, since the centralizer of a maximal torus is $T$ itself.  
It is also expressed in the form $G_{\C}/P$, where $G_{\C}$
denotes the complexification of $G$ and $P$ a parabolic subgroup of $G_{\C}$.  Furthermore, it is realized 
as an orbit of the coadjoint representation of $G$ on $\mathfrak{g}^{*}$, the dual of the Lie algebra $\mathfrak{g}$
of $G$.  From these realizations,  a generalized flag manifold $G/H$  has many interesting geometric properties. 
In fact, $G/H$ admits a complex structure, a K\"{a}hler structure and a symplectic structure, 
and therefore plays an important role in algebraic topology, differential geometry,  
and algebraic geometry.

In algebraic topology,  it is a classical problem to determine the integral cohomology ring $H^*(G/H;\Z)$ of 
a generalized flag manifold $G/H$, and many mathematicians tried to find the method to compute it.  
Here, we briefly review the history of determination of the cohomology ring 
of $G/H$; The first notable  result is due to  Borel. In his thesis \cite{Bor53}, 
Borel  described the rational cohomology ring $H^*(G/H;\Q)$ of $G/H$ (as of all homogeneous spaces $G/H$, with 
$H$ of maximal rank in $G$) in terms of the rings of invariants of the Weyl groups of $G$ and $H$ (see also 
\S \ref{rational}). This description of $H^*(G/H;\Q)$ is called the {\it Borel presentation}.  
 However, Borel's result does not hold in general when $G$ or $H$ has torsion. 
  Using the Morse theory,  Bott showed that $G/H$ is free of torsion, i.e., $H^{*}(G/H;\Z)$ is a free $\Z$-module, 
 and its odd Betti numbers vanish (\cite[Theorem A]{Bott56}).\footnote{This fact also follows from the cell decomposition
  of $G/H \cong G_{\C}/P$ by {\it Schubert cells} (\cite{Che58}).}  Furthermore, there is an algorithm for constructing the ring $H^*(G/H;\Z)$ in terms of the Cartan integers  (\cite[Theorem III'']{Bott-Sam58}). 
So the problem of computing  the integral cohomology ring of a generalized flag manifold $G/H$ could be  solved in principle. 
However, in practice, it becomes unmanageable to carry out the above algorithm.  
After that, another new technique was introduced by   Baum \cite{Bau68}. 
 He used the Eilenberg-Moore spectral sequence converging to $H^*(G/H;k)$, where $k$ is a field or the integers. 
 Under certain hypothesis, he showed that the Eilenberg-Moore spectral sequence of the fibration $G/H \longrightarrow 
 BH \longrightarrow BG$ collapses, where $BH$ (resp. $BG$) denotes the classifying space of $H$ (resp. $G$). 
 However, this collapse theorem does not give the ring structure of $H^*(G/H;k)$ unless we solve  the 
 so-called {\it extension problem} of a spectral sequence (see \cite[Chapter 8]{McC01}).   
 Although these general methods are quite useful when we discuss general properties of the cohomology ring
 of $G/H$, it seems  difficult to apply them directly to some homogeneous spaces of exceptional Lie groups, 
  especially of $E_{8}$.   Using the Borel presentation of the rational cohomology ring $H^{*}(G/H;\Q)$ and 
  the results on mod $p$ cohomology rings $H^*(G;\Z/p\Z)$ for all primes $p$,  Toda  initiated the research  for 
  computing the integral cohomology ring of $G/H$ with $H$ torsion free (\cite{Toda75}). 
    Along the line of  his idea, the integral cohomology rings of various flag manifolds have been computed 
    explicitly 
   (\cite{Toda-Wat74}, \cite{Wat75}, \cite{Ishi-Toda77}, \cite{Ishi77}, \cite{Wat98}, 
    \cite{Nak01}).\footnote{Some of these results contain minor errors that are corrected in \cite{DZ05}, \cite{Kaji-Nak1}.}  
    Toda's method is quite useful in practical computations. For instance, based on Toda's result, 
    Kono and  Ishitoya computed the mod $2$  cohomology ring $H^*(E_{8}/T;\Z/2\Z)$ explicitly 
    (\cite{Kono-Ishi87}),\footnote{Their result also contains 
    an error that was corrected in \cite{Tot05}.} and Totaro computed the $\Z_{(2)}$-cohomology ring 
    $H^*(E_{8}/A_{8};\Z_{(2)})$, where $A_{8}$ denotes the subgroup of $E_{8}$ isomorphic to $SU(9)/(\Z/3\Z)$. 
    Using these  results, Totaro  succeeded in computing the {\it torsion index}\footnote{The torsion index is an 
    integer associated to any compact connected Lie group defined by  Grothendieck \cite{Gro58}.} of $E_{8}$ (\cite{Tot05}).   
    Finally, we should mention another direction towards determination of $H^*(G/H;\Z)$; 
    As mentioned before, there is a description of $H^*(G/H;\Z)$ in terms of the cell decomposition of $G/H \cong G_{\C}/P$ 
   by  Schubert cells.
    This originated in the work of Ehresmann \cite{Ehr34}  on the cell decomposition of a complex Grassman manifold, 
    and was later extended to general flag manifolds by  Chevalley \cite{Che58}.  
    This description of $H^*(G/H;\Z)$ is called the {\it Schubert presentation}.  In this description, 
     the additive structure of $H^*(G/H;\Z)$  is described in terms of the so-called  {\it Schubert classes}
     indexed by certain subset of the Weyl group of $G$ (see \cite[\S 5]{BGG73}). 
     As for the multiplicative structure, we have to compute the  {\it structure constants} for the multiplication 
     of Schubert classes (see for example \cite{Pra05}).   
     This is one of the main problems of the {\it Schubert calculus}.\footnote{In the case of a complex Grassman manifold, 
      the multiplicative rule of Schubert classes is  known as   the {\it Littlewood-Richardson rule} (see \cite{Ful97}).}   
     In this direction,  many authors studied the Chow rings of generalized flag manifolds 
    (\cite{CMP08}, \cite{Ili-Man05}, \cite{NS06},  \cite{NSZ09}).\footnote{The  Chow ring  $A(G/H)$  of $G/H$  is canonically   
    isomorphic to $H^*(G/H;\Z)$ via the cycle map (\cite[\S 6]{Gro58}), \cite[Chapter 19]{Ful98}).} 
    Recently  Duan developed extensively a multiplicative rule of   Schubert
    classes which is a generalization of the Littlewood-Richardson rule of a complex Grassmann manifold (\cite{Duan05}). 
    Furthermore,  he and  Zhao  computed the integral cohomology rings of the above flag manifolds 
    independently of Toda's method (\cite{DZ05}, \cite{DZ06}, \cite{DZ08}).\footnote{The relations between their results and ours 
     are revealed  in \cite{Kaji-Nak1} via the {\it divided difference operators} due to Bernstein-Gelfand-Gelfand  \cite{BGG73}
     and Demazure \cite{Dem73}.}

  Until recently  none of these  methods have been successful in computing 
  the {\it integral}  cohomology rings of homogeneous spaces of the exceptional 
  Lie group $E_8$.\footnote{Partial result was obtained by Nikolenko and Semenov  \cite{NS06}.} 
    The group $E_8$   contains a closed connected subgroup $T^1 \! \cdot \! E_7$  whose local type is $T^1 \times E_7$, 
     where $T^1$  is a  one dimensional torus (see \cite[\S 2]{Ishi-Toda77}).  It is obtained as 
     the centralizer of a certain one dimensional torus (see \ref{notations}). Hence the homogeneous space 
     $E_8/T^1 \! \cdot \! E_7$ is a generalized flag manifold. 
     In this paper, using the above method due to Borel and Toda, we compute the integral cohomology ring of 
      $E_8/T^1 \! \cdot \! E_7$ explicitly.\footnote{Duan and Zhao also computed $H^*(E_{8}/T^1 \! \cdot \! E_{7};\Z)$ 
      in terms of Schubert classes (\cite[Theorem 7]{DZ05}.}

    The motivation of the current work is not only  the determination of  the integral cohomology ring 
   $H^*(E_{8}/T^1 \! \cdot \! E_{7};\Z)$ itself, but  the cohomology of the irreducible symmetric space 
   $E IX = E_8/S^3 \! \cdot \! E_7$ (see \cite[TABLE V]{Hel01}, \cite[\S 2]{Ishi-Toda77}), as well as the integral cohomology ring 
   of the full flag manifold $E_8/T$.  The symmetric space $E IX$ is a {\it quaternionic K\"{a}hler manifold}
  and $E_{8}/T^1 \! \cdot \! E_{7}$ is the {\it twister space} of $E IX$  (\cite{Sal82}). 
  Using the Gysin exact sequence 
  associated with the $2$-sphere bundle $S^{2} \longrightarrow E_{8}/T^1 \! \cdot \! E_{7} \longrightarrow E IX$, 
  we can compute the integral cohomology ring of $E IX$.\footnote{The real cohomology rings
  of quaternionic K\"{a}hler manifolds are studied in \cite{Nag-Tak87}.}  Using  Theorem \ref{thm:E_8/T^1E_7} of this paper, 
  we also compute   $H^{*}(E_8/T;\Z)$,  which is the only remaining case
    among $G/T$'s  for $G$ simple, in our forthcoming paper (\cite{Nak2}). 
   In \cite{Kaji-Nak2}, the result of $H^*(E_{8}/T;\Z)$ will be  used to compute  the Chow ring $A(\mathrm{E}_{8})$  of 
   the complex algebraic  group $\mathrm{E}_{8}$
   (For the Chow rings of complex algebraic groups, see \cite{Gro58}, \cite{Kac85}, \cite{Mar74}).     
   Moreover,   the homogeneous space $E_8/T^1 \! \cdot \! E_7$ is a {\it generating variety}  of $E_8$ in the sense of 
   Bott \cite{Bott58} and its integral cohomology ring is  needed to compute the Pontrjagin ring 
   $H_{*}(\Omega E_8;\Z)$, where $\Omega E_8$ denotes the based loop  space of $E_8$.  
 
 The paper is  organized as follows:   In \S2,  we compute the rings of invariants of 
 the Weyl groups of $E_8$ and $ T^1 \! \cdot \! E_7$ which are needed for the computations of the rational cohomology ring of 
 $E_{8}/T^1 \! \cdot \! E_{7}$.
    In \S3, using these results and the Borel presentation
 of a rational cohomology ring, we compute the rational cohomology ring of $E_8/T^1 \! \cdot \! E_7$. 
 In \S 4,   by  investigating the integral cohomology ring of $E_8/T$ in low degrees, we determine 
 the integral cohomology ring  of $E_8/T^1 \! \cdot \! E_7$ explicitly (Theorem \ref{thm:E_8/T^1E_7}). 
 Furthermore, in \S 5,  using the Gysin 
 exact sequence associated with the  circle bundle $S^1 \longrightarrow E_8/E_7  \longrightarrow 
 E_8/T^1 \! \cdot \! E_7$,  we also  determine the integral cohomology  of $E_8/E_7$ (Corollary  \ref{cor:E_8/E_7}).

 {\it Acknowledgments}. We thank  Mamoru Mimura for many useful suggestions and advice
 during the preparation of the manuscript. Thanks are also due to  Haibao Duan and Xhezi Zhao,
 who informed us their result on the Schubert presentation of $H^*(E_{8}/T^1 \! \cdot \! E_{7};\Z)$. 
 Finally, we thank Tetsu Nishimoto and Shizuo Kaji for valuable discussions.  
 The main part of this manuscript was obtained around 2007. I would like to apologize for not posting it sooner.

\section{Rings of invariants of Weyl groups }
In this section, we compute the rings of invariants of the Weyl groups of $E_{8}$ and 
$T^1 \! \cdot \! E_{7}$ over $\Q$ explicitly that are needed for the computation of the rational 
cohomology ring of $E_{8}/T^1 \! \cdot \! E_{7}$.  
\subsection{Notations}  \label{notations}
Let $T$ be a fixed maximal torus of $E_{8}$. According to \cite{Bou68},  the Dynkin diagram of $E_{8}$ is as follows:

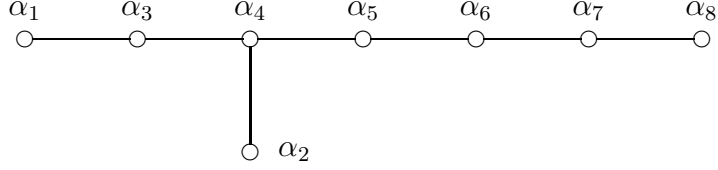
\setlength{\unitlength}{1mm}
     \begin{picture}(100,30)
        \multiput(30,20)(15,0){7}{\circle{2}}     
        \put(60,5){\circle{2}}     
        \multiput(31,20)(15,0){6}{\line(1,0){13}}
        \put(60,19){\line(0,-1){13}}          
        \put(30,25){\makebox(0,0)[t]{$\alpha_{1}$}}    
        \put(45,25){\makebox(0,0)[t]{$\alpha_{3}$}}
        \put(60,25){\makebox(0,0)[t]{$\alpha_{4}$}}
        \put(75,25){\makebox(0,0)[t]{$\alpha_{5}$}}
        \put(90,25){\makebox(0,0)[t]{$\alpha_{6}$}}
        \put(105,25){\makebox(0,0)[t]{$\alpha_{7}$}}
        \put(120,25){\makebox(0,0)[t]{$\alpha_{8}$}}
         \put(68,5){\makebox(0,0)[r]{$\alpha_{2}$}}    
   \end{picture} \\
where $\alpha_{i} \; (1 \leq i \leq 8)$ are the simple roots.\footnote{In topology, a {\it root}  $\alpha$ is a real linear form on 
$\mathfrak{t}$, the Lie algebra of $T$ (\cite[Chapter 4]{Ada69}. It is customary that $2 \pi \sqrt{-1} \alpha$, regarded as 
a complex linear form on $\mathfrak{t}_{\C}$, the complexification of $\mathfrak{t}$, is 
a {\it root} in the theory of complex Lie algebras.}  As usual we may regard each
root as an element of $H^{1}(T;\Z)$\footnote{Let $\exp: \mathfrak{t} \longrightarrow T$ be the exponetial map. The {\it integrer lattice} (or {\it unit lattice})  $\Gamma$ is defined as the inverse image of the unit $e$ of $T$ under $\exp$ (\cite[4.11]{Ada69}). 
A real linear form $\alpha$ on $\mathfrak{t}$ is said to be  {\it integral} if it takes integral values on $\Gamma$ 
(\cite[1.2]{Bor-Hir58}).    Then we may identify  the integral real linear forms on $\mathfrak{t}$  with
 $\Hom (\Gamma, \Z) \cong H^{1}(T;\Z)$.  Notice that each root takes integral values on $\Gamma$ by definition.}, and therefore an element of $H^{2}(BT;\Z)$ via the negative {\it transgression} (see \cite[\S 10]{Bor-Hir58}).

Let $C_{8}$ be the centralizer of a one dimensional torus determined by
$\alpha_{i} = 0 \; (i \neq 8)$ in $\mathfrak{t}$. Then as shown in  \cite[\S 2]{Ishi-Toda77}, 
    \[  C_{8} = T^1 \! \cdot \! E_7  \quad \text{and} \quad T^1 \cap E_7 \cong \Z/2\Z.   \] 
The Weyl groups  of $E_{8}$ and $C_{8}$ are respectively given  as follows:
   \[   W(E_{8}) = \langle s_{i} \; (1 \leq i \leq 8) \rangle, \quad 
        W(C_{8}) = \langle s_{i} \; (i \leq i \leq 7)  \rangle,    \]
where $s_{i}$ denotes the simple reflection corresponding to the simple root $\alpha_{i}$.  

Let $\{ \omega_{i} \}_{1 \leq i \leq 8}$ be the fundamental weights corresponding 
to the system of the simple roots $\{ \alpha_{i} \}_{1 \leq i \leq 8}$. 
We also regard each weight as an elememt of $H^{2}(BT;\Z)$. Then 
$\{ \omega_{i} \}_{1 \leq i \leq 8}$ forms a basis of $H^{2}(BT;\Z)$ and 
$H^{*}(BT;\Z) = \Z[\omega_{1}, \omega_{2}, \ldots, \omega_{8}]$. The action 
of $s_{i}$'s on $\{ \omega_{i} \}_{1 \leq i \leq 8}$ is given  as follows:
 \[   s_{i}(\omega_{k}) =  \left \{ 
                \begin{array}{llll} 
                           & \hspace{-0.3cm}   \omega_{i} - \displaystyle{\sum_{j=1}^{8} \frac{2 (\alpha_{i}|\alpha_{j})}
                          {(\alpha_{j}|\alpha_{j})} \omega_{j} } \quad & \text{if}  \quad k = i, \\ 
                           & \hspace{-0.3cm}  \omega_{k}  \quad                         & \text{if} \quad k \neq i,  
                \end{array}   \right.   \]
where $(\; \cdot \; | \; \cdot \;)$ denotes the $W(E_{8})$-invariant inner product on $\mathfrak{t}^{*}$, the dual of $\mathfrak{t}$.

Now  we introduce  the elements of $H^{2}(BT;\Z)$ as follows
 (Throughout this paper,  $\sigma_{i}(x_{1}, \ldots, x_{n})$ denotes the $i$-th elementary 
  symmetric function in the variables $x_1, \ldots, x_{n}$):  
  \begin{align*}
    t_{8} &= \omega_{8}, \\
    t_{7} &= s_{8}(t_{8})  = \omega_{7} - \omega_{8}, \\
    t_{6} &= s_{7}(t_{7})  = \omega_{6} - \omega_{7}, \\
    t_{5} &= s_{6}(t_{6})  = \omega_{5} - \omega_{6}, \\
    t_{4} &= s_{5}(t_{5})  = \omega_{4} - \omega_{5}, \\
    t_{3} &= s_{4}(t_{4})  = \omega_{2} + \omega_{3} - \omega_{4}, \\
    t_{2} &= s_{3}(t_{3})  = \omega_{1} + \omega_{2} - \omega_{3}, \\
    t_{1} &= s_{1}(t_{2})  = -\omega_{1} + \omega_{2},  
  \end{align*} 
  \begin{align*} 
    c_{i} &= \sigma_{i}(t_{1}, \ldots,  t_{8}), \\ 
    t     &= \omega_{2}  = \frac{1}{3}c_{1}.
   \end{align*}
Then $t$ and $t_{i}$'s span $H^{2}(BT;\Z)$,  since each $\omega_{i}$ is an integral 
linear combination of $t$ and $t_{i}$'s,  and we have the following  isomorphism:
  \[  H^{*}(BT;\Z) = \Z [t_{1},\ldots,t_{8},t]/(c_{1}-3t).   \]
The elements $\{ t_{i} \}_{1 \leq i \leq 8}$ and $t$ behave nicely under the action of 
the Weyl group $W(E_{8})$ of $E_{8}$. In fact, the action of $s_{i} \; (1 \leq i \leq 8)$  on $\{ t_{i} \}_{1 \leq i \leq 8}$ and $t$ is 
 given by T{\scriptsize ABLE} 1, where blanks indicate the trivial action.

\begin{table}[htbp] 
   \begin{center}  
   \begin{tabular}{|c|c|c|c|c|c|c|c|c|}
   \noalign{\hrule height0.8pt}
   \hfil $ $ &  $s_{1}$ & $s_{2} $ & $ s_{3}$ & $s_{4}$ & $s_{5}$ &  $s_{6}$ & $s_{7}$ & $s_{8}$   \\
   \hline
  $t_{1}$ & $t_{2}$ &  $ t-t_{2}-t_{3}$ & $$ & $$ & $$ & $$ & $$ & $$    \\  
  \hline
 $t_{2}$ &  $t_{1}$ &  $t-t_{1}-t_{3}$ & $t_{3}$ & $$ & $$ & $$ & $$ & $$  \\  
  \hline
   $t_{3}$ & $$       &  $t-t_{1}-t_{2}$ & $t_{2}$ & $t_{4}$  & $$ & $$ &  $$ & $$   \\
  \hline     
   $t_{4}$ & $$ & $$ & $$  & $t_{3}$ & $t_{5}$ & $$ & $$ & $$   \\
  \hline
   $t_{5}$ & $$ & $$ & $$ & $$ & $t_{4}$ & $t_{6}$ & $$ & $$   \\
  \hline
   $t_{6}$ & $$ & $$ & $$ & $$ & $$ & $t_{5}$ & $t_{7}$ & $$   \\
  \hline
   $t_{7}$ & $$ & $$ & $$ & $$ & $$ & $$ & $t_{6}$ & $t_{8}$   \\
   \hline
   $t_{8}$ & $$ & $$ & $$ & $$ & $$ & $$ & $$ & $t_{7}$   \\  
    \hline
   $t$     & $$ & $ 2t-t_{1}-t_{2}-t_{3}$ & $$ & $$ & $$ & $$ & $$ & $$   \\
   \noalign{\hrule height0.8pt}
   \end{tabular} 
   \end{center}
   \caption{  }
   \end{table}

Next we introduce a  basis of $H^{2}(BT;\Q)$ which behaves nicely under the action of the 
Weyl group $W(C_{8})$ of $C_{8}$: 
  \begin{align*}
      u        &= t_{8},  \\ 
     \tau      &= t - \frac{3}{2}u, \\ 
     \tau_{i}  & = t_{i}- \frac{1}{2}u \; (1 \leq i \leq 7).   
  \end{align*} 
Then we have 
  \[ H^{*}(BT;\Q) = \Q[u , \tau, \tau_{1}, \ldots,  \tau_{7}]/( \bar{c}_{1}-3\tau)  
                  = \Q[u,  \tau_{1}, \ldots, \tau_{7}]  \]
for $\bar{c}_{1} = \tau_{1} + \cdots + \tau_{7}$.   The action of $s_{i} \; (1 \leq i \leq 7)$  on 
$\{ \tau_{i} \}_{1 \leq i \leq 7}$ and $\tau$ is  given by T{\scriptsize ABLE} 2,  where  blanks
also  indicate the trivial action.

   \begin{table}[htbp]
   \begin{center}  
   \begin{tabular}{|c|c|c|c|c|c|c|c|}
   \noalign{\hrule height0.8pt}
   \hfil $ $ &  $s_{1}$ & $s_{2} $ & $ s_{3}$ & $s_{4}$ & $s_{5}$ & $s_{6}$ & $s_{7}$   \\
   \hline
   $\tau_{1}$ & $\tau_{2}$ &  $ \tau-\tau_{2}-\tau_{3}$ & $$ & $$ & $$ & $$ & $$   \\
   \hline  
   $\tau_{2}$ &  $\tau_{1}$ &  $\tau-\tau_{1}-\tau_{3}$ & $\tau_{3}$ & $$ & $$  & $$ & $$    \\
   \hline  
   $\tau_{3}$ & $$       &  $\tau-\tau_{1}-\tau_{2}$ & $\tau_{2}$ & $\tau_{4}$  & $$ & $$ &  $$   \\
   \hline     
   $\tau_{4}$ & $$ & $$ & $$  & $\tau_{3}$ & $\tau_{5}$ & $$ & $$  \\
    \hline
   $\tau_{5}$ & $$ & $$ & $$ & $$ & $\tau_{4}$ & $\tau_{6}$ & $$  \\
   \hline
   $\tau_{6}$ & $$ & $$ & $$ & $$ & $$ & $\tau_{5}$ & $\tau_{7}$   \\
   \hline
   $\tau_{7}$ & $$ & $$ & $$ & $$ & $$ & $$ & $\tau_{6}$ \\
   \hline
   $\tau$     & $$ & $ -\tau+\tau_{4}+\tau_{5}+\tau_{6}+\tau_{7}$ & $$ & $$ & $$ & $$ & $$ \\
   \hline
   $u$ & $$ & $$ & $$ & $$ & $$ & $$  & $ $   \\
   \noalign{\hrule height0.8pt}
   \end{tabular} 
   \end{center}
   \caption{ }
   \end{table}

Since $E_{7} \cap T = T'$ is a maximal torus of $E_{7}$, we have a 
commutative diagram of natural maps:
\begin{equation}
\begin{CD}
   E_{7}/T' @>{\sim}>> C_{8}/T @>{i}>> E_{8}/T   \\
   @VVV
    @VVV
   @VV{\iota_{0}}V   \\
    BT' @>{g}>> BT @>{ =}>> BT.
\end{CD}
\end{equation}
Since $E_{8}$ is 2-connected, $\iota_{0}^{*}: H^{2}(BT;\Z)  \longrightarrow   H^{2}(E_{8}/T;\Z)$ 
is an isomorphism.  Under this isomorphism, we denote the $\iota_{0}^{*}$-images of 
$t_{i} \; (1 \leq i \leq 8)$ and  $t$  by the same letters. Thus we have the generators 
$t_{i} \; (1 \leq i \leq 8)$ and $t$ of $H^{2}(E_{8}/T;\Z)$ with a relation $c_{1} = 3t$. 
We donote the generators $t_{i} \; (1 \leq i \leq 7)$  and $x$ in \cite[\S 1]{Wat75} by $t_{i}' \;(1 \leq i \leq 7)$ and $t'$
respectively. 
Then,   by a  similar argument to that  in \cite[\S 1]{Wat75}, we have
 \begin{equation} \label{eqn:t_i}  
 g^{*}(t_{i}) = t_{i}' \; (1 \leq i \leq 7), \; g^{*}(t_{8}) = 0, \;
     g^{*}(t) = t'
  \end{equation}  
under the homomorphism $g^{*}: H^{2}(BT;\Z) \longrightarrow H^{2}(BT';\Z)$.

\subsection{Ring of invariants of  $W(C_8)$}   \label{W(C_8)-invariants} 
In this subsection, we compute the ring of invariants of the Weyl group $W(C_{8})$ explicitly. 
First we recall the rational invariant forms for $E_{7}$ given in \cite[\S 2]{Wat75} (see also \cite[2.2]{Meh88}). 
We put  
  \[ x_{i}' = 2t_{i}'-t' \; (1 \leq i \leq 7)  \quad \text{and} \quad  x_{8}' = t'.  \]
Then the set 
  \[ S' = \{ \pm \; (x_{i}'+x_{j}') \; ( 1 \leq i < j \leq 8 )  \}   \subset H^{2}(BT' ; \Q)  \]
is invariant under the action of $W(E_{7})$. Thus we have $W(E_{7})$-invariant forms
 \begin{equation}  \label{inv.forms.E7}
     I_{n}' = \sum_{y \in S'}y^{n}  \in H^{2n}(BT';\Q)^{W(E_{7})} \quad (n \geq 0).  
 \end{equation}
Then direct computation using the method as in  \cite[\S 2]{Wat75}  yields the
following results:
\begin{equation} \label{eqn:inv.forms.E_7}   
\begin{array}{cl} 
   I_{2}'    &= -2^5 \cdot 3 (c_{2}' - 4{t'}^{2}), \medskip   \\ 
   I_{6}'    &\equiv 2^8 \cdot 3^2 ( {c_{3}'}^{2} + 8c_{6}') 
                      \mod (t', \frak{a}_{6}' ),  \medskip \\
    I_{8}'   &\equiv 2^{12 }\cdot 5 (2{c_{4}'}^{2} - 3c_{3}'c_{5}') 
                      \mod (t', \frak{a}_{8}'),  \medskip \\
    I_{10}'  &\equiv 2^{12} \cdot 3^2 \cdot 5 \cdot 7 ({c_{5}'}^{2}-
                     4c_{3}'c_{7}')  \mod (t', \frak{a}_{10}'), \medskip  \\
    I_{12}'  &\equiv 2^{15} \cdot 3^2 \cdot 5 (-54{c_{6}'}^{2} 
                     + 18c_{5}'c_{7}' -c_{3}'c_{4}'c_{5}')  \mod (t', \frak{a}_{12}'),  \medskip  \\
    I_{14}'  &\equiv 2^{16} \cdot 3 \cdot 7 \cdot 11 \cdot 29 (2{c_{7}'}^{2}
                     +2c_{3}'c_{4}'c_{7}' - c_{3}'c_{5}'c_{6}') \mod (t', \frak{a}_{14}'), \medskip \\
    I_{18}'  &\equiv 2^{21} \cdot 5 \cdot 1229 (-126c_{5}'c_{6}'c_{7}'
                     -5c_{3}'c_{4}'c_{5}'c_{6}') \mod (t', \frak{a}_{18}'), 
  \end{array}
\end{equation}
where $c_{i}' = \sigma_{i}(t_{1}', \ldots ,t_{7}')$ and $\frak{a}_{i}'$
denotes the ideal of $H^{*}(BT';\Q)$ generated by $I_{j}'$ for $j <i$ with  $j \in \{ 2,6,8,10,12,14,18 \}$.
We also recall  the following result:
 \begin{prop}[\cite{Wat75},  Lemma 2.1, \cite{Meh88}, 2.2]   \label{prop:W(E_7)-invariants}    
  The ring of invariants of the Weyl group $W(E_{7})$  over $\Q$  is given  as follows$:$
   \[ H^{*}(BT';\Q)^{W(E_{7})} 
    = \Q[I_{2}' ,I_{6}',I_{8}',I_{10}',I_{12}',I_{14}',I_{18}'].  \]
 \end{prop}  

 T{\scriptsize ABLE} 2 shows that the action of $W(C_{8})$ on $\tau,\tau_{1},\ldots,\tau_{7}$ is the same as 
  that of $W(E_{7})$ on $t',t_{1}',\ldots,t_{7}'$.  Therefore if we put 
  \[  \chi_{i} = 2\tau_{i} - \tau \; (1 \leq i \leq 7)  \quad \text{and} \quad \chi_{8} = \tau,  \] 
  the set 
  \[  \Sigma = \{ \pm (\chi_{i} + \chi_{j}) \; (1 \leq i <  j \leq 8) \} \subset H^{2}(BT;\Q) 
  \] 
  is invariant under the action of $W(C_{8})$.  We define $W(C_{8})$-invariant forms $J_{n} \; (n \geq 0)$ 
  as 
  \[    J_{n} = \sum_{y \in \Sigma} y^{n} \in H^{2n}(BT;\Q)^{W(C_{8})}.   
  \]
  Then, by Proposition \ref{prop:W(E_7)-invariants}, we have 
  \begin{lem} \label{lem:W(C_8)-invariants}
   The ring of invariants  of the Weyl group $W(C_{8})$ over $\Q$ is  given as follows$:$
    \[  H^{*}(BT;\Q)^{W(C_{8})} =  
                                \Q[u, J_{2},J_{6}, J_{8},J_{10},J_{12},J_{14},J_{18}].   
   \]
 \end{lem}

\subsection{Ring of invariants of  $W(E_8)$}
In this subsection, we compute the ring of invariants of the Weyl group $W(E_{8})$ over $\Q$  explicitly. 
According to Chevalley  \cite{Che55}, this ring of invariants is generated by $8$ algebraically 
independent polynomials (basic invariants) $f_{1}, \ldots, f_{8}$  of degrees $2, 8, 12, 14, 18, 20, 24, 30$ 
(\cite[p.59]{Hum90}).    We will give these invariant polynomials explicitly.  
 
We put 
   \[  \xi_{i} = 2t_{i} - \frac{2}{3}t \; (1 \leq i \leq 8) \quad \text{and} \quad  
       \xi_{9} = -\frac{2}{3}t.      
   \] 
Then the action of $s_{i} \; (1 \leq i \leq 8)$ on $\{  \xi_{i} \}_{1 \leq i \leq 9}$ is given by 
T{\scriptsize ABLE} 3, where $\eta = (\xi_{1} + \xi_{2} + \xi_{3})/3$. 

\begin{table}[htbp] 
   \begin{center}  
   \begin{tabular}{|c|c|c|c|c|c|c|c|c|}
   \noalign{\hrule height0.8pt}
   \hfil $ $  &  $s_{1}$    &  $s_{2}$            & $s_{3}$   & $s_{4}$   & $s_{5}$ &  $s_{6}$ & $s_{7}$ & $s_{8}$   \\
   \hline
   $\xi_{1}$  &  $\xi_{2}$  &  $\xi_{1} - 2\eta$  & $$        & $$        & $$ & $$ & $$ & $$    \\  
   \hline
   $\xi_{2}$  &  $\xi_{1}$  &  $\xi_{2} - 2\eta$  & $\xi_{3}$ & $$        & $$ & $$ & $$ & $$  \\  
   \hline
   $\xi_{3}$  &  $$         &  $\xi_{3} - 2\eta$  & $\xi_{2}$ & $\xi_{4}$ & $$ & $$ &  $$ & $$   \\
   \hline     
   $\xi_{4}$  &  $$         &  $\xi_{4} + \eta$   & $$        & $\xi_{3}$ & $\xi_{5}$ & $$ & $$ & $$   \\
   \hline
   $\xi_{5}$  &  $$         &  $\xi_{5} + \eta$   & $$        & $$        & $\xi_{4}$ & $t_{6}$ & $$ & $$   \\
   \hline
   $\xi_{6}$  &  $$         &  $\xi_{6} + \eta$   & $$        & $$        & $$ & $\xi_{5}$ & $\xi_{7}$ & $$   \\
   \hline
   $\xi_{7}$  &  $$         &  $\xi_{7} + \eta$   & $$        & $$        & $$ & $$ & $\xi_{6}$ & $\xi_{8}$   \\
   \hline
   $\xi_{8}$  &  $$         &  $\xi_{8} + \eta$   & $$        & $$        & $$ & $$ & $$ & $\xi_{7}$   \\  
    \hline
   $\xi_{9}$  &  $$         &  $\xi_{9} + \eta$   & $$        & $$        & $$ & $$ & $$ & $$   \\
   \noalign{\hrule height0.8pt}
   \end{tabular} 
   \end{center}
   \caption{  }
   \end{table} 
   
From T{\scriptsize ABLE} 3, we see  that the set
  \[ S = \{ \pm \;(\xi_{i}-\xi_{j}) \; (1 \leq i < j \leq 9) , \; 
            \pm \; (\xi_{i}+\xi_{j}+\xi_{k}) \; (1 \leq i < j < k \leq 9) \}  \subset H^{2}(BT;\Q)   \] 
is invariant under the action of $W(E_{8})$.\footnote{In fact, $S$ is an orbit of $2 \omega_{8}$ under the action of 
$W(E_{8})$.  Since $\omega_{8}$ is equal to the highest root $\tilde{\alpha}$ (\cite{Bou68}), it turns out that 
$S$ is $2$ times  the root system of $E_{8}$.}    Thus we have $W(E_{8})$-invariant forms\footnote{For a compact 
connected Lie group $G$, the set of roots is nothing but the set of non-zero weights of the adjoint representation of $G$ on $\mathfrak{g}$. 
Therefore  $I_{n}$ is equal to the $n$-th 
Chern character of the adjoint representation of $E_{8}$ up to the multiple of $2$.}
  \[  I_{n} = \sum_{ y \in S}y^{n}  \in H^{2n}(BT;\Q)^{W(E_{8})} \quad (n \geq 0).   \]

Let us compute $I_{n}$'s in the following way;   We put
   \[ s_{n} = \xi_{1}^{n} + \cdots + \xi_{9}^{n}, \quad 
      d_{n} = \sigma_{n}(\xi_{1},\ldots,\xi_{9}).    \]
Then $s_{n}$'s and $d_{n}$'s are related to each other by the Newton formula:
  \begin{equation} \label{eqn:Newton}
      s_{n} = \sum_{i=1}^{n-1} (-1)^{i-1}s_{n-i}d_{i}+(-1)^{n-1}nd_{n}
        \hspace{0.2cm}  (d_{n} = 0 \hspace{0.2cm}  \text{for} \hspace{0.2cm} n > 9 ). 
  \end{equation} 
Note that $s_{0} = 9$ and 
   \[ s_{1} = d_{1} = \sum_{i=1}^{9} \xi_{i} = \sum_{i=1}^{8} \left (2t_{i} - \dfrac{2}{3}t \right ) - \dfrac{2}{3} t = 2(c_{1}  - 3t) =  0.  \]
Then 
  \begin{align*}
      \sum_{n \geq 0} \frac{I_{n}}{n!}  
           & =  \sum_{1 \leq i < j \leq 9}e^{\xi_{i}-\xi_{j}} + \sum_{1 \leq i < j \leq 9}e^{-\xi_{i}+\xi_{j}} 
            + \sum_{1 \leq i < j < k \leq 9}e^{\xi_{i}+\xi_{j}+\xi_{k}}  
            + \sum_{1 \leq i < j < k \leq 9}e^{-\xi_{i}-\xi_{j}-\xi_{k}} \\
           & =    \left  (\sum_{i = 1}^{9} e^{\xi_{i}} \right )  \left (\sum_{j = 1}^{9} e^{-\xi_{j}} \right )  - 9  
               +\frac{1}{3} \left( \sum_{i = 1}^{9} e^{3\xi_{i}} + \sum_{i = 1}^{9} e^{-3\xi_{i}} \right ) \\ 
           & - \frac{1}{2}   \left \{ \left  (\sum_{i = 1}^{9} e^{\xi_{i}} \right  ) 
               \left  (\sum_{i = 1}^{9} e^{2\xi_{i}} \right )  + \left (\sum_{i = 1}^{9} e^{-\xi_{i}} \right)
               \left (\sum_{i = 1}^{9} e^{-2\xi_{i}}\right   ) \right  \}  \\
           & +\frac{1}{6} \left \{ \left (\sum_{i = 1}^{9} e^{\xi_{i}} \right )^{3} 
             + \left (\sum_{i = 1}^{9} e^{-\xi_{i}} \right )^{3}\right  \}  \\
           & =  \left (\sum_{n \geq 0}\frac{s_{n}}{n!} \right )  \left (\sum_{m \geq 0}\frac{(-1)^{m}s_{m}}{m!} \right )-9
            + \frac{1}{3} \left ( \sum_{n \geq 0}\frac{3^{n}s_{n}}{n!} 
            + \sum_{n \geq 0}\frac{(-1)^{n}3^{n}s_{n}}{n!} \right ) \\
           & -\frac{1}{2} \left \{ \left (\sum_{n \geq 0}\frac{s_{n}}{n!} \right )
             \left (\sum_{m \geq 0}\frac{2^{m}s_{m}}{m!} \right )
            + \left (\sum_{n \geq 0}\frac{(-1)^{n}s_{n}}{n!} \right )
              \left (\sum_{m \geq 0}\frac{(-1)^{m}2^{m}s_{m}}{m!} \right ) \right \}  \\
            & +\frac{1}{6}  \left \{ \left (\sum_{n \geq 0}\frac{s_{n}}{n!} \right )^{3} 
               + \left (\sum_{n \geq 0}\frac{(-1)^{n}s_{n}}{n!} \right )^{3}  \right   \}. 
\end{align*}                               
Therefore we have
\begin{equation}  \label{eqn:inv.forms.E_8}
 \begin{array}{cll} 
    I_{n}  &  =   \displaystyle{\sum_{i = 0}^{n} \binom{n}{i}(-1)^{n-i} s_{i}s_{n-i} 
                + 2 \cdot 3^{n-1} s_{n}  }  
                - \displaystyle{\sum_{i = 0}^{n} \binom{n}{i}2^{n-i}s_{i}s_{n-i} }  \medskip \\ 
           &   + \dfrac{1}{3} \displaystyle{\sum_{i = 0}^{n} \sum_{j = 0}^{n-i} \binom{n}{i} 
   \binom{n-i}{j} s_{i}s_{j}s_{n-i-j}  }
 \end{array}
\end{equation} 
for $n$ even.   On the other hand,  since $\xi_{i} = 2t_{i}-\dfrac{2}{3}t \; (1 \leq i \leq 8)$
and $ \xi_{9} = -\dfrac{2}{3}t$,  we have
  \begin{align*} 
      \sum_{n=0}^{9}d_{n} &= \prod_{i=1}^{9}(1+\xi_{i}) 
             = \left (1-\frac{2}{3}t \right ) \prod_{i=1}^{8}
              \left (1-\frac{2}{3}t+2t_{i} \right ) \\
     & =  \left (1-\frac{2}{3}t \right ) \sum_{i=0}^{8}
         \left (1-\frac{2}{3}t \right )^{8-i}2^{i}c_{i}
  \end{align*}
and therefore 
    \begin{equation}  \label{eqn:d_n} 
         d_{n} = \sum_{i=0}^{n}\left \{ \binom{8-i}{n-i}
        +\binom{8-i}{n-i-1} \right \}
          \left (-\frac{2}{3}t \right )^{n-i}2^{i}c_{i}   
     \end{equation} 
 for $1 \leq n \leq 9$.  Using (\ref{eqn:inv.forms.E_8}), (\ref{eqn:Newton}) and (\ref{eqn:d_n}), 
 $I_{n}$ can be expressed as a polynomial in $t$ and $c_{i} \; (2 \leq i \leq 8)$, 
 and we obtain  the following results:
\begin{lem} \label{lem:inv.forms.E_8}  
 \begin{equation*}
  \begin{array}{cll}  
   \mathrm{(i)} \;      I_{2}  & =       -2^5 \cdot 3 \cdot 5  (c_{2}-4t^{2}),  \medskip \\
   \mathrm{(ii)} \;     I_{8}  & \equiv   2^{14} \cdot 3 \cdot 5  \{ -18c_{8}-3c_{3}c_{5}+2c_{4}^{2} + t(12c_{7}-3c_{3}c_{4})  
                                 + t^2(-6c_{6} +3c_{3}^{2})  \medskip \\
                               & + 12t^{3}c_{5}+2t^{4}c_{4}-12t^{5}c_{3}+14t^{8}\} 
                                  \mod   \tilde{\frak{a}}_{8},  \medskip \\
   \mathrm{(iii)}\;     I_{12}  & \equiv    2^{18} \cdot 3^{5} \cdot 7  
                                    \left (c_{6}^{2}-\dfrac{5}{3}c_{5}c_{7} +\dfrac{5}{54}c_{3}c_{4}c_{5}
                                   -\dfrac{1}{6}c_{3}^{2}c_{6}+ \dfrac{1}{24}c_{3}^{4} \right )  
                                   \mod (t,c_{8},  \tilde{\frak{a}}_{12}),   \medskip  \\
   \mathrm{(iv)} \;     I_{14}  & \equiv   2^{20} \cdot 3^{2} \cdot 5^{2} \cdot 7 
                                 \cdot 11  \left (c_{7}^{2}  -\dfrac{1}{2}c_{3}c_{5}c_{6}
                                 + \dfrac{1}{3}c_{3}c_{4}c_{7}+ \dfrac{1}{6}c_{4}c_{5}^{2} \right )  
                                 \mod (t,c_{8},  \tilde{\frak{a}}_{14}),  \medskip   \\
   \mathrm{(v)} \;       I_{18}  & \equiv   - 2^{23} \cdot 3^{2} \cdot 5 \cdot 7 \cdot 13 
                                 \left ( c_{3}^{6} -7c_{3}^{4}c_{6} +\dfrac{29}{9}c_{3}^{3}c_{4}c_{5}
                                              +182c_{3}^{2}c_{5}c_{7}  +75c_{3}c_{5}^{3}  \right.  \medskip \\
                                           &  \left.  -\dfrac{476}{3}c_{3}c_{4}c_{5}c_{6}
                                              -24c_{5}c_{6}c_{7}   
                                          \right )  \mod (t,c_{8}, \tilde{\frak{a}}_{18}),  \medskip  \\
   \mathrm{(vi)} \;     I_{20} & \equiv  2^{27} \cdot 3^{3} \cdot 5^{2} \cdot 11   \cdot 17 \cdot 41
                                \left ( \dfrac{1}{144}c_{5}^{4}- \dfrac{1}{18}c_{3}c_{5}^{2}c_{7}
                               -\dfrac{1}{54}c_{3}^{2}c_{4}c_{5}^{2}   -\dfrac{1}{27}c_{3}^{3}c_{4}c_{7}  \right. \medskip \\
                              &   \left. + \dfrac{1}{18}c_{3}^{3}c_{5}c_{6} \right ) 
                                  \mod (t,c_{8} , \tilde{\frak{a}}_{20}),  \medskip   \\
   \mathrm{(vii)}\; I_{24}  & \equiv   2^{32} \cdot 3^3 \cdot 5^2 \cdot 7 \cdot 
                          11 \cdot 19  \cdot 199 \left ( \dfrac{31}{8640}c_{3}^5c_{4}c_{5}
                          + \dfrac{1}{480}c_{3}^{4}c_{5}c_{7}  +\dfrac{337}{25920}c_{3}^{3}c_{5}^{3}  \right. \medskip \\
                          & \left.  - \dfrac{71}{4320}c_{3}^{3}c_{4}c_{5}c_{6} 
                          + \dfrac{31}{240}c_{3}^{2}c_{5}c_{6}c_{7} 
                          +\dfrac{31}{480}c_{3}c_{5}^{3}c_{6}   -\dfrac{22}{135}c_{3}c_{4}c_{5}^{2}c_{7} 
                          - \dfrac{1}{120}c_{4}c_{5}^{4} \right )  \medskip \\
                          & \mod (t, c_{8},   \tilde{\frak{a}}_{20}),   \medskip        \\
   \mathrm{(viii)} \;  I_{30} & \equiv   2^{38} \cdot 3^4 \cdot 5^5 \cdot 7 \cdot 
                          11 \cdot 13 \cdot 61  \left ( -\dfrac{599}{51840}c_{3}^{5}c_{4}c_{5}c_{6}
                          + \dfrac{47}{34560}c_{3}^{5}c_{5}^{3}   \right.  \medskip \\ 
                          & + \dfrac{1519}{25920}c_{3}^{4}c_{5}c_{6}c_{7}  
                          + \dfrac{6293}{7290}c_{3}^{3}c_{4}c_{5}^{2}c_{7} 
                          - \dfrac{32537}{25920}c_{3}^{3}c_{5}^{3}c_{6}
                          + \dfrac{189919}{466560}c_{3}^{2}c_{4}c_{5}^{4}   \medskip  \\
                          &   \left.  +\dfrac{2012}{1215}c_{3}c_{4}c_{5}^{2}c_{6}c_{7} 
                          -\dfrac{16693}{25920}c_{3}c_{5}^{4}c_{7} 
                          -\dfrac{223}{6480}c_{4}c_{5}^{4}c_{6} 
                          -\dfrac{1}{1728}c_{5}^{6} \right )   \mod (t, c_{8},  \tilde{\frak{a}}_{20}),   \medskip  
 \end{array}
\end{equation*}
where $c_{i} = \sigma_{i}(t_{1}, \ldots ,t_{8})$ and $\tilde{\frak{a}}_{i}$
denotes the ideal of $H^{*}(BT;\Q)$ generated by $I_{j}$ for $j <i$ with  $j \in \{ 2,8,12,14,18, 20, 24, 30 \}$.
\end{lem}

Next consider the following elements of $H^{*}(BT;\Q)^{W(C_{8})} \subset H^{*}(BT;\Q)$:
\begin{equation} \label{eqn:relations}
\begin{array}{cl}
    \tilde{I}_{20} & =    9u^{20} + 45u^{14}v + 12u^{10}w + 60u^8v^2 + 30u^4vw   + 10u^2v^3 + 3w^2, \medskip  \\
    \tilde{I}_{24} & =  11u^{24} + 60u^{18}v +  21u^{14}w + 105u^{12}v^2 + 60u^8vw   + 60u^{6}v^3  + 9u^4 w^2 
                       \medskip \\
                  &  + 30 u^2 v^2 w + 5v^4, \medskip   \\
    \tilde{I}_{30} & =   -9u^{30} - 24u^{24}v  -12u^{20}w  + 36u^{14}vw -40u^{12}v^3  -12u^{10}w^2  \medskip \\
                  &   + 120u^8 v^2w   -140u^{6}v^{4} + 24u^{4}vw^{2} -40u^2v^3w   -16v^5 -8w^{3},   \medskip   
   \end{array}
\end{equation}
where 
\begin{equation} \label{eqn:generators}
\begin{array}{cl}  
  u &= t_{8}, \medskip  \\
  v &= \dfrac{1}{46080} J_{6} -\dfrac{273}{640} u^{6}, \medskip  \\ 
  w &= \dfrac{1}{15482880} J_{10}
      -\dfrac{55}{24} u^{4}v -\dfrac{666919}{645120} u^{10}.  \medskip 
\end{array}
\end{equation}

\begin{rem} 
  By the classical result of Borel \cite{Bor53}, the rational cohomology ring of $E_{8}/C_{8}$ is 
  described in terms of the rings of invariants of the Weyl groups$:$  
  \[     H^*(E_{8}/C_{8};\Q) \cong H^*(BT;\Q)^{W(C_{8})}/(H^{+}(BT;\Q)^{W(E_{8})})   
  \]
   $($see $\S \ref{rational})$. Under this isomorphism, the above elements $u$, $v$, $w$ represent the rational 
    cohomology classes in $H^{*}(E_{8}/C_{8};\Q)$.        
     Moreover,   the integral cohomology ring of $E_{8}/C_{8}$ is  torsion free $($\cite{Bott56}$)$, 
     and hence, is contained in  the rational cohomology ring$:$ $H^*(E_{8}/C_{8};\Z) \hookrightarrow H^*(E_{8}/C_{8};\Q)$. 
    Then the elements $u$, $v$, $w$ are in fact the integral cohomology classes in $H^{*}(E_{8}/C_{8};\Z)$ 
   $($For details, see $\ref{u,v,w,x})$.  
\end{rem}

We wish to find  the relations among $W(C_{8})$-invariants $J_{n}$, $W(E_{8})$-invariants $I_{n}$  and $\tilde{I}_{20}$,
 $\tilde{I}_{24}$, $\tilde{I}_{30}$.  For this purpose, we consider the ring of invariants: 
\[  A = H^{*}(BT;\Q)^{ \langle s_{1},s_{3},\ldots, s_{7}  \rangle }. 
\] 
Then $A$ is a subalgebra of $H^{*}(BT;\Q)$ containing both $H^{*}(BT;\Q)^{W(C_{8})}$ and $H^*(BT;\Q)^{W(E_{8})}$.  
More explicitly,   we have
  \begin{equation} \label{eqn:A}
     A  = \Q[u, t, c_{2}, \ldots, c_{7}].
  \end{equation} 
In fact,  we can show (\ref{eqn:A}) as follows;   Putting  $\tilde{c}_{i} = \sigma_{i}(t_{1},\ldots,t_{7})$, 
we have 
  \[ c_{n} = \tilde{c}_{n} + u\tilde{c}_{n-1}  \;  (1 \leq n \leq 8),  \]
since 
 \[ \sum_{n = 0}^{8} c_{n} = \prod_{i = 1}^{8}(1 + t_{i}) 
  = (1 + u) \prod_{i = 1}^{7}(1 + t_{i}) 
  = (1 + u) \sum_{n = 0}^{7}\tilde{c}_{n}.   \]  
Conversely,  one can express 
 \[  \tilde{c}_{n} = c_{n} - uc_{n-1} + u^{2}c_{n-2} - \cdots +(-1)^{n}u^{n}  \; (1 \leq n \leq 7).  \]
In particular, the following relation holds:
 \begin{equation} \label{eqn:c_{8}}
     c_{8} = uc_{7}-u^{2}c_{6}+u^{3}c_{5}- u^{4}c_{4} + u^5c_{3} - u^6 c_{2} + u^7 c_{1} -u^{8}. 
 \end{equation}
Therefore, by T{\scriptsize ABLE} 1, we have  
 \begin{align*}
     A &= H^{*}(BT;\Q)^{ \langle s_{1},s_{3},\ldots,s_{7} \rangle} \\
       & = \Q[t_{1},t_{2},\ldots,t_{7},u]^{ \langle s_{1},s_{3}, \ldots,s_{7} \rangle}  \\
       & = \Q[u,\tilde{c}_{1},\tilde{c}_{2},\ldots,\tilde{c}_{7}] \\
       & =\Q[u,c_{1},c_{2},\ldots,c_{7}]   \\
       & = \Q[u, t, c_{2}, \ldots, c_{7}],  
 \end{align*}
which has shown (\ref{eqn:A}). 

Denote by 
  \[  \mathfrak{a}_{i} \subset A \; (\text{resp.}  \; \mathfrak{b}_{i}   \subset  
      H^{*}(BT;\Q)^{W(C_{8})}),   \]
the ideal of $A$  (resp. of $H^{*}(BT;\Q)^{W(C_{8})}$)  generated
by $I_{j}$'s for $j < i$ where $j \in \{ 2,8,12,14,18,20,24,30 \}$.

The remainder  of this section is devoted to proving  the next lemma:  
\begin{lem} \label{lem:invariants}
   In $H^{*}(BT;\Q)^{W(C_{8})} = \Q[u,J_{2},J_{6},J_{8},J_{10}, J_{12},J_{14},J_{18}]$,  we have
 \begin{enumerate} 
 \item [(i)] \begin{align*} 
        & I_{2} = 5J_{2} + 120 u^2,   \\   
        & I_{8} = 2^2 \cdot 3J_{8} + (\text{decomp.}),   \\ 
        & I_{12} = -2^{2} \cdot 7 J_{12} + (\text{decomp.}),   \\
        & I_{14} = \dfrac{2^{3} \cdot 3 \cdot 5^{2}}{29}J_{14} + (\text{decomp.}),   \\
        & I_{18} = -\dfrac{2^{4} \cdot 3 \cdot 5^{2} \cdot 7 \cdot 13}{1229} J_{18} + (\text{decomp.}),      
    \end{align*} 
  where $($decomp.$)$ means decomposable elements of  $\Q[u, J_{2}, J_{6}, J_{8}, J_{10}, J_{12}, J_{14}, J_{18}]$. 
 \item [(ii)] 
      \begin{align*}    
         & I_{20}  \equiv 2^{27} \cdot 3^2 \cdot 5^2 \cdot 11
                   \cdot 17   \cdot 41 \tilde{I}_{20}   \mod \mathfrak{b}_{20},   \\  
         & I_{24}  \equiv 2^{32} \cdot 3^3 \cdot 5 \cdot 7 \cdot 11
                   \cdot 19 \cdot 199  \tilde{I}_{24}  \mod \mathfrak{b}_{20},  \\
         & I_{30}  \equiv 2^{35} \cdot 3^{4} \cdot 5^{5} \cdot 7 \cdot  11 \cdot 13 \cdot 61 
                   \tilde{I}_{30} \mod \mathfrak{b}_{20}. 
    \end{align*} 
 \end{enumerate} 
\end{lem}

Suppose Lemma \ref{lem:invariants} for the moment. Then we have the following (see also \cite[2.3]{Meh88}):
\begin{lem}  \label{lem:W(E_8)-invariants}
    The ring of invariants  of the Weyl group $W(E_{8})$ over $\Q$ is 
given as follows$:$
\[  H^{*}(BT;\Q)^{W(E_{8})} = \Q[I_{2},I_{8},I_{12}, I_{14},I_{18},I_{20},I_{24},I_{30}].    
\]
\end{lem}

\begin{proof}
  By Lemma \ref{lem:invariants},  $I_{i} \; (i = 2, 8, 12, 14, 18, 20, 24, 30)$  are algebraically independent. 
  On the other hand, by the result of Borel \cite{Bor53},   
  the ring of invariants $H^*(BT;\Q)^{W(E_{8})}$ is isomorphic to  the rational cohomology ring $H^*(BE_{8};\Q)$ 
  of the classifying space $BE_{8}$, and it is known that  
\[    H^{*}(BE_{8};\Q) \cong   \Q[y_{4},y_{16},y_{24},y_{28},y_{36}, y_{40},y_{48},y_{60}]     
\]
with deg ($y_{i}) = i$.   Therefore we have the required result.
\end{proof}

\begin{proof}[Proof of Lemma $\ref{lem:invariants} \;  \mathrm{(i)}$]
Since $\tau_{i} = t_{i} - \dfrac{1}{2}u \; (1 \leq i \leq 7)$ and  
$\tau = t - \dfrac{3}{2}u $, we have 
\[   \tau_{i} \equiv t_{i} \;(1 \leq i \leq 7), \; \tau \equiv t \mod(u). 
\] 
Therefore, putting  $\bar{c}_{i} = \sigma_{i}(\tau_{1},\ldots ,\tau_{7}) \; (1 \leq i \leq 7)$,  
we obtain  
\begin{equation}  \label{eqn:bar{c}_n}   
  c_{n} \equiv \bar{c}_{n} \; (1 \leq n \leq 7), \; c_{8} \equiv 0 \mod (u).   
\end{equation} 
Since $J_{n}$ has the same expression as $I_{n}'$ by simply replacing $t'$, $c_{i}'$ with $\tau$, $\bar{c}_{i}$ 
(see \ref{W(C_8)-invariants}),  we have, by (\ref{eqn:inv.forms.E_7}) and (\ref{eqn:bar{c}_n}), 
\begin{equation} \label{eqn:J_i}
  \begin{array}{cl}
     J_{2} & =  -2^5 \cdot 3   (\bar{c}_{2} - 4\tau^{2})  \medskip \\ 
           & \equiv -2^5 \cdot 3   c_{2}   \mod (t,u),   \medskip \\
    J_{6}  & \equiv 2^8 \cdot 3^2 ({\bar{c}_{3}}^{2} + 8 \bar{c}_{6})  
           \mod (\tau, \bar{\frak{a}}_{6}) \medskip  \\ 
           &   \equiv 2^8 \cdot 3^2   (c_{3}^{2} + 8c_{6}) \mod (t,u,\bar{\frak{a}}_{6}), \medskip  \\
    J_{8}  & \equiv 2^{12} \cdot 5 (2{\bar{c}_{4}}^{2} - 3\bar{c}_{3}\bar{c}_{5}) 
             \mod (\tau, \bar{\frak{a}}_{8}) \medskip  \\
           & \equiv 2^{12} \cdot 5  (2c_{4}^{2} - 3c_{3}c_{5}) \mod (t,u,\bar{\frak{a}}_{8}), \medskip \\
    J_{10} & \equiv 2^{12} \cdot 3^2 \cdot 5 \cdot 7 
          ({\bar{c}_{5}}^{2}-4\bar{c}_{3}\bar{c}_{7}) 
          \mod (\tau,\bar{\frak{a}}_{10}) \medskip \\
        & \equiv 2^{12} \cdot 3^2 \cdot 5 \cdot 7 (c_{5}^{2} - 4c_{3}c_{7})
                         \mod(t,u, \bar{\frak{a}}_{10}),  \medskip     \\
    J_{12} &\equiv 2^{15} \cdot 3^2 \cdot 5   (-54{\bar{c}_{6}}^{2} +18 \bar{c}_{5}\bar{c}_{7} 
                -\bar{c}_{3}\bar{c}_{4}\bar{c}_{5})   
                  \mod (\tau, \bar{\frak{a}}_{12})  \medskip  \\
            &\equiv 2^{15} \cdot 3^2 \cdot 5 (-54c_{6}^{2}+18c_{5}c_{7}-c_{3}c_{4}c_{5})
                  \mod (t,u,\bar{\frak{a}}_{12}),  \medskip \\
    J_{14} &\equiv 2^{16} \cdot 3 \cdot 7 \cdot 11 \cdot 29 (2{\bar{c}_{7}}^{2}
            + 2\bar{c}_{3}\bar{c}_{4}\bar{c}_{7} 
                             - \bar{c}_{3}\bar{c}_{5}\bar{c}_{6})
                      \mod (\tau, \bar{\frak{a}}_{14}) \medskip  \\ 
        &\equiv 2^{16} \cdot 3 \cdot 7 \cdot 11 \cdot 29 
                       (2c_{7}^{2}+2c_{3}c_{4}c_{7} - c_{3}c_{5}c_{6})
                      \mod (t,u, \bar{\frak{a}}_{14}),  \medskip \\
    J_{18} &\equiv 2^{21} \cdot 5 \cdot 1229 
         (-126\bar{c}_{5}\bar{c}_{6}\bar{c}_{7} - 
             5\bar{c}_{3}\bar{c}_{4}\bar{c}_{5}\bar{c}_{6})  
                 \mod (\tau,\bar{\frak{a}}_{18})  \medskip \\
        &  \equiv 2^{21} \cdot 5 \cdot 1229 
         (-126c_{5}c_{6}c_{7} -5c_{3}c_{4}c_{5}c_{6})
                   \mod (t,u, \bar{\frak{a}}_{18}), 
\end{array}
\end{equation}
where $\bar{\frak{a}}_{i}$ denotes the ideal of $H^{*}(BT;\Q)^{W(C_{8})}$ 
generated by $J_{j}$'s for $j < i$ with  $j \in \{2,6,8,10,12,14,18 \}$.

Now we prove the last formula of (i).  
Since $I_{18} \in H^*(BT;\Q)^{W(E_{8})} \subset H^{*}(BT;\Q)^{W(C_{8})} \\ 
 = \Q[u,J_{2},J_{6},J_{8},J_{10}, J_{12},J_{14},J_{18}]$,   
we can put 
  \begin{equation} \label{eqn:I_18}
   I_{18} =  \alpha_{18} J_{18}  + (\text{decomp.})   
  \end{equation}
 for some $\alpha_{18} \in \Q$. On the other hand,  by using 
(\ref{eqn:A}) and  (\ref{eqn:J_i}),  we have
  \begin{align*}
     A/(t,u,\bar{\frak{a}}_{18}) = & \; A/(t,u,J_{2},J_{6},J_{8},J_{10},
   J_{12},J_{14}) \\
      = & \; \Q[u, t, c_{2}, c_{3}, c_{4}, c_{5}, c_{6}, c_{7}] \\
        & \hspace{0.5cm} 
            \left /\left (  \begin{array}{lllll} 
                       & \hspace{-0.4cm} t, \; u, \;  c_{2}, \; c_{3}^{2}+8c_{6}, \; c_{4}^{2} - \dfrac{3}{2}c_{3}c_{5}, \;  
                          c_{5}^{2}-4c_{3}c_{7}, \\
                       & \hspace{-0.4cm} c_{6}^{2} - \dfrac{1}{3}c_{5}c_{7} + \dfrac{1}{54}c_{3}c_{4}c_{5}, \;  
                         c_{7}^{2} + c_{3}c_{4}c_{7} - \dfrac{1}{2}c_{3}c_{5}c_{6} 
                       \end{array}  \right )  \right.  \\
      = &  \; \Q[c_{3}, c_{4}, c_{5}, c_{6}, c_{7}] \\
        & \hspace{0.5cm} 
          \left / \left ( \begin{array}{llll} 
                          & \hspace{-0.4cm}   c_{3}^{2}+8c_{6}, \; 
                           c_{4}^{2} - \dfrac{3}{2}c_{3}c_{5}, \;  c_{5}^{2}-4c_{3}c_{7},  \\
                          & \hspace{-0.4cm} c_{6}^{2} - \dfrac{1}{3}c_{5}c_{7} + \dfrac{1}{54}c_{3}c_{4}c_{5}, \;
                           c_{7}^{2} + c_{3}c_{4}c_{7} - \dfrac{1}{2}c_{3}c_{5}c_{6} 
                          \end{array}  \right ) \right. .
  \end{align*}   
We consider (\ref{eqn:I_18}) in the ring $A/(t,u,\bar{\frak{a}}_{18})$. 
 Then,  by Lemma \ref{lem:inv.forms.E_8} (v) and (\ref{eqn:J_i}), we have   
  \begin{align*}
      I_{18} &\equiv -2^{25} \cdot 3 \cdot 5^{3} \cdot 7 \cdot 13 
           (-126c_{5}c_{6}c_{7} -5c_{3}c_{4}c_{5}c_{6} ), \\
       J_{18} &\equiv 2^{21} \cdot 5 \cdot 1229 
              (-126c_{5}c_{6}c_{7} -5c_{3}c_{4}c_{5}c_{6}).       
  \end{align*}
Therefore $\alpha_{18} = -  2^{4} \cdot 3 \cdot 5^{2} \cdot 7 \cdot 13/1229$. 
A similar tedious computation gives the other formulas.
\end{proof}

Before proceeding the proof of Lemma \ref{lem:invariants} (ii),  we need the 
following lemma:
\begin{lem} \label{lem:J6,J10}
  Explicit forms of $W(C_{8})$-invariants $J_{6},J_{10}$ are given as follows$:$
  \begin{align*}
      J_{6} & =  2^{8} \cdot 3 \{24 \bar{c}_{6} + 3 {\bar{c}_{3}}^{2} 
              -4 \bar{c}_{2}\bar{c}_{4} -2\bar{c}_{2}^{3} +(-12\bar{c}_{5} 
             -6 \bar{c}_{2}\bar{c}_{3}) \tau +(31{\bar{c}_{2}}^{2} + 16\bar{c}_{4}) {\tau}^{2} \\
      &      +12\bar{c}_{3}{\tau}^{3} -136 \bar{c}_{2} {\tau}^{4} + 188 {\tau}^{6} \},  \\
      J_{10} & =  2^{12} \cdot 3 \{ 105 {\bar{c}_{5}}^{2} 
       -420 \bar{c}_{3}\bar{c}_{7} + 90 \bar{c}_{2}\bar{c}_{3}\bar{c}_{5}
       -60 \bar{c}_{2}{\bar{c}_{4}}^{2} + 300 {\bar{c}_{2}}^{2}\bar{c}_{6}    
       +15 {\bar{c}_{2}}^{2}{\bar{c}_{3}}^{2} \\
     & -20 {\bar{c}_{2}}^{3}\bar{c}_{4} -2 {\bar{c}_{2}}^{5} 
       + (-30 {\bar{c}_{2}}^{3}\bar{c}_{3} - 330 {\bar{c}_{2}}^{2}\bar{c}_{5}
       +480 \bar{c}_{2}\bar{c}_{7} + 90 \bar{c}_{2}\bar{c}_{3}\bar{c}_{4}
       -210 \bar{c}_{4}\bar{c}_{5} )\tau  \\
    & +(270 {\bar{c}_{2}}^{2}\bar{c}_{4} -210 \bar{c}_{2}{\bar{c}_{3}}^{2}
      -150 \bar{c}_{3}\bar{c}_{5} -2220 \bar{c}_{2}\bar{c}_{6}
      +75 {\bar{c}_{2}}^{4} + 345 {\bar{c}_{4}}^{2} ) {\tau}^{2} \\
    & +(480 {\bar{c_{2}}}^{2}\bar{c_{3}} -570 \bar{c_{3}}\bar{c_{4}} 
      +2070 \bar{c}_{2}\bar{c}_{5} -660 \bar{c}_{7} ) {\tau}^{3}  \\
    & +(-1050 \bar{c}_{2}\bar{c}_{4} + 4080 \bar{c}_{6} 
      -950 {\bar{c}_{2}}^{3} +705 {\bar{c}_{3}}^{2} ) {\tau}^{4} 
   +(-2250 \bar{c}_{2}\bar{c}_{3} -3420 \bar{c}_{5} ) {\tau}^{5}  \\
   &  +(1580 \bar{c}_{4} + 5165 {\bar{c}_{2}}^{2}) {\tau}^{6} +2820 \bar{c}_{3} {\tau}^{7} 
    -12360 \bar{c}_{2} {\tau}^{8} + 10868 {\tau}^{10} \},  
\end{align*}   
where   $\bar{c}_{i} = \sigma_{i}( \tau_{1},  \ldots , \tau_{7} )$.
\end{lem}
\begin{proof}
   Since $J_{n}$ has the same expression as $I_{n}'$ by replacing $t',c_{i}'$ with
$\tau, {\bar{c}}_{i}$ (see \ref{W(C_8)-invariants}), we have to compute
the $W(E_{7})$-invariant forms $I_{6}',I_{10}'$ explicitly. But this can be
done from the data in   \cite[\S 2]{Wat75}.
\end{proof}

In order to show Lemma \ref{lem:invariants} (ii), we need to describe the elements $v$, $w$ (see (\ref{eqn:generators})) 
in the ring $A$, in other words, in terms of $t$, $u$, $c_{i} \; (2 \leq i \leq 7)$. 
First will  rewrite $J_{6}$, $J_{10}$ in terms of $t$, $u$, $c_{i} \; (2 \leq i \leq 7)$; 
Since $u = t_{8}$ and $\tau_{i} = t_{i} - \dfrac{1}{2}u \; (1 \leq i \leq 7)$,
 we have 
\begin{align*}
   \left ( 1 + \frac{1}{2}u \right ) \sum_{n = 0}^{7} \bar{c}_{n} 
    &=   \left (1 + \frac{1}{2}u \right ) \prod_{i = 1}^{7}(1 + \tau_{i})  
    =   \left (1 + \frac{1}{2}u \right )\prod_{i = 1}^{7}
        \left (1 -\frac{1}{2}u + t_{i} \right )  \\
     &=  \prod_{i = 1}^{8} \left (1 - \frac{1}{2}u + t_{i} \right ) 
      = \sum_{i = 0}^{8} \left (1 - \frac{1}{2}u \right )^{8-i} c_{i},   
\end{align*}
and hence 
   \[    \bar{c}_{n} + \frac{1}{2}u {\bar{c}_{n-1}} 
       =  \sum_{i = 0}^{n} 
         \binom{8-i}{n-i} \left (-\frac{1}{2} \right )^{n-i} c_{i} u^{n-i} \quad 
         (1 \leq n \leq 7).   \]
From this,   we obtain
\begin{equation}   \label{eqn:bar{c_i}}
  \begin{array}{cll} 
     \bar{c}_{1} &= 3t - \dfrac{9}{2}u,  \medskip   \\
    \bar{c}_{2} &= c_{2} - 12tu +\dfrac{37}{4}u^{2},   \medskip \\
    \bar{c}_{3} &= c_{3} -\dfrac{7}{2}c_{2}u + \dfrac{87}{4}tu^{2}
                   -\dfrac{93}{8}u^{3}, \medskip  \\
    \bar{c}_{4} &= c_{4} -3c_{3}u +\dfrac{11}{2}c_{2}u^{2} -24 tu^{3}
                   + \dfrac{163}{16}u^{4}, \medskip  \\
    \bar{c}_{5} &= c_{5} - \dfrac{5}{2}c_{4}u + 4c_{3}u^{2} 
            -\dfrac{21}{4}c_{2}u^{3} + \dfrac{297}{16}tu^{4} - 
             \dfrac{219}{32}u^{5},  \medskip \\
    \bar{c}_{6} &= c_{6} - 2c_{5}u + \dfrac{11}{4}c_{4}u^{2} 
            - \dfrac{13}{4}c_{3}u^{3} + \dfrac{57}{16}c_{2}u^{4}
            - \dfrac{45}{4}tu^{5} + \dfrac{247}{64}u^{6},  \medskip \\
    \bar{c}_{7} &= c_{7} - \dfrac{3}{2}c_{6}u + \dfrac{7}{4}c_{5}u^2- \dfrac{15}{8}c_{4}u^{3}
            + \dfrac{31}{16}c_{3}u^{4} - \dfrac{63}{32}c_{2}u^{5}
            + \dfrac{381}{64}tu^{6}   \medskip \\
             & - \dfrac{255}{128}u^{7}.
  \end{array} 
\end{equation} 
By Lemma \ref{lem:J6,J10} and (\ref{eqn:bar{c_i}}), we can express $J_{6}$, $J_{10}$ (and 
hence $v$, $w$) in terms of $t$, $u$, $c_{i} \; (2 \leq i \leq 7)$.  In fact, we need only 
the expressions modulo certain elements (see (\ref{eqn:v.w}) below);  
By Lemma \ref{lem:inv.forms.E_8} (ii)  and (\ref{eqn:c_{8}}),   we have 
   \begin{align*}
      I_{8} &\equiv 2^{14} \cdot 3 \cdot 
   5(2c_{4}^{2}-3c_{3}c_{5}) \mod (t,c_{8},\frak{a}_{8}),  \\
      c_{8} &= uc_{7}-u^{2}c_{6}+u^{3}c_{5}-u^{4}c_{4}
          +u^{5}c_{3}-u^{6}c_{2}+u^{7}c_{1}-u^{8} \\
            &\equiv uc_{7}-u^{2}c_{6}+u^{3}c_{5}-u^{4}c_{4}+u^{5}c_{3}-u^{8}
              \mod (t,\frak{a}_{8}), 
     \end{align*}
and hence
\begin{equation}  \label{eqn:mod(t.c_8.I_2.I_8)}    
      \begin{array}{llll} 
    c_{4}^{2} &\equiv \dfrac{3}{2}c_{3}c_{5} \mod (t,c_{8},\frak{a}_{12}),  \medskip \\
    u^{8} &\equiv uc_{7}-u^{2}c_{6}+u^{3}c_{5}-u^{4}c_{4}+u^{5}c_{3}
               \mod (t,c_{8},\frak{a}_{8}). \medskip  
    \end{array}
\end{equation}  
Therefore,  by  (\ref{eqn:generators}), Lemma \ref{lem:J6,J10}, (\ref{eqn:bar{c_i}}) and (\ref{eqn:mod(t.c_8.I_2.I_8)}), 
 we obtain
\begin{equation}  \label{eqn:v.w}
  \begin{array}{clll}
     v &= \dfrac{1}{46080}J_{6} - \dfrac{273}{640}u^{6}  \medskip \\
       & \equiv \dfrac{2}{5}c_{6} + \dfrac{1}{20}c_{3}^{2} - \dfrac{1}{2}c_{5}u
         + \dfrac{1}{3}c_{4}u^{2} - \dfrac{1}{2}c_{3}u^{3}  \mod (t , \frak{a}_{8}), \medskip  \\  
      w &= \dfrac{1}{15482880}J_{10} -\dfrac{55}{24}u^{4}v 
           -\dfrac{666919}{645120}u^{10} \medskip  \\
        & \equiv \dfrac{1}{12}c_{5}^{2} - \dfrac{1}{3}c_{3}c_{7} 
          + \left (\dfrac{1}{2}c_{3}c_{6}-\dfrac{1}{6}c_{4}c_{5} \right )u
          -\dfrac{1}{6}c_{3}c_{5}u^{2}
          + \left (-c_{7}+\dfrac{1}{3}c_{3}c_{4} \right )u^{3}  \medskip \\
       &  -\dfrac{1}{2}c_{3}^{2}u^{4}  + \dfrac{1}{3}c_{4}u^{6} 
          + \dfrac{1}{2}c_{3}u^{7}   \mod (t, c_{8}, \frak{a}_{12}).
 \end{array}                         
\end{equation} 
Under these preparations,  we will prove Lemma \ref{lem:invariants} (ii). 

\begin{proof}[Proof of  Lemma $\ref{lem:invariants}  \; \mathrm{(ii)}$]
First note that 
   \begin{align*}
       H^{*}(BT;\Q)^{W(C_{8})} & = \Q[u,J_{2},J_{6}, J_{8},J_{10},J_{12},J_{14},J_{18}]  \\
                               & = \Q[u,I_{2},v,I_{8},w,I_{12}, I_{14},I_{18}]
   \end{align*}  by  (i) and (\ref{eqn:generators}). 
Since $I_{20} \in H^{*}(BT;\Q)^{W(E_{8})} \subset H^{*}(BT;\Q)^{W(C_{8})}$, we can put
 \begin{align*}
  (*) \quad  I_{20} & \equiv \;  2^{27}\cdot 3^{3} \cdot 5^{2}  \cdot 11 \cdot 17 \cdot 41
            (\lambda_{1}u^{20}+\lambda_{2}u^{14}v+\lambda_{3}u^{8}v^{2}  + \lambda_{4}u^{2}v^{3}  \\
            & + \lambda_{5}u^{10}w + \lambda_{6}u^{4}vw +  \lambda_{7}w^{2} ) \mod \frak{b}_{20}            
 \end{align*}
for some $\lambda_{i} \in \Q$. 
In order to determine the coefficients $\lambda_{i}$, we need  the following lemma, 
which is directly verified by making use of (\ref{eqn:A}) and  Lemma \ref{lem:inv.forms.E_8}: 
\begin{lem} \label{lem:kakikae}
 \begin{align*} 
           A/(t,c_{8},\frak{a}_{20}) & = \;  A/(t,I_{2},c_{8},I_{8},I_{12},I_{14},I_{18})  \\
                                     & = \;  \Q[u,c_{3},c_{4},c_{5},c_{6},c_{7}]/J,  
 \end{align*} 
         where $J$ is the ideal generated by 
             \[  
                      \begin{array}{llll} 
                      & u^{8} - uc_{7} + u^{2}c_{6} - u^{3}c_{5} + u^{4}c_{4} - u^{5}c_{3},  \medskip \\
                      & c_{4}^{2}-\dfrac{3}{2}c_{3}c_{5},   \medskip \\
                      & c_{6}^{2}-\dfrac{5}{3}c_{5}c_{7}
                      +\dfrac{5}{54}c_{3}c_{4}c_{5}-\dfrac{1}{6}c_{3}^{2}c_{6}
                      +\dfrac{1}{24}c_{3}^{4}, \medskip \\
                      &  c_{7}^{2}-\dfrac{1}{2}c_{3}c_{5}c_{6}+\dfrac{1}{3}c_{3}c_{4}c_{7}
                      +\dfrac{1}{6}c_{4}c_{5}^{2},  \medskip \\
                      & c_{3}^{6}-7c_{3}^{4}c_{6}+\dfrac{29}{9}c_{3}^{3}c_{4}c_{5}
                        +182c_{3}^{2}c_{5}c_{7} + 75c_{3}c_{5}^{3} 
                        -\dfrac{476}{3}c_{3}c_{4}c_{5}c_{6}-24c_{5}c_{6}c_{7}.  \medskip 
                   \end{array} 
             \] 
In particular,   $A/(t,c_{8},\frak{a}_{20})$ has a basis
  $\{ u^{i}c_{3}^{j}c_{4}^{k}c_{5}^{l}c_{6}^{m}c_{7}^{n} 
   \; (0 \leq i \leq 7,\; 0 \leq j \leq 5, \; 0 \leq l,\; 0 \leq k,m,n \leq 1 )\}$ 
as a $\Q$-vector space.
\end{lem}

Now we consider the relation $(*)$ in the ring $A/(t,c_{8},\frak{a}_{20})$. By Lemma \ref{lem:inv.forms.E_8} 
 (vi), we have
  \[  I_{20} \equiv 2^{27} \cdot 3^{3} \cdot 5^{2}  \cdot 11 \cdot 17 \cdot 41
       \left    (\frac{1}{144}c_{5}^{4}-\frac{1}{18}c_{3}c_{5}^{2}c_{7}
      -\frac{1}{54}c_{3}^{2}c_{4}c_{5}^{2}-\frac{1}{27}c_{3}^{3}c_{4}c_{7}
      +\frac{1}{18}c_{3}^{3}c_{5}c_{6} \right).    \]
On the other hand,  using (\ref{eqn:v.w}) and Lemma \ref{lem:kakikae}, we can rewrite each  monomial
in the right hand side  of $(*)$. For example,  we have 
  \begin{align*}
    w^{2}   &\equiv   \frac{1}{144}c_{5}^{4}-\frac{1}{18}c_{3}c_{5}^{2}c_{7}
                -\frac{1}{54}c_{3}^{2}c_{4}c_{5}^{2}
                -\frac{1}{27}c_{3}^{3}c_{4}c_{7} +\frac{1}{18}c_{3}^{3}c_{5}c_{6} -\frac{1}{12}u^{7}c_{3}^{2}c_{7}  \\
            & +\frac{1}{12}u^{7}c_{3}c_{5}^{2}+
              \frac{1}{24}u^{6}c_{3}^{3}c_{5}+\frac{1}{3}u^{6}c_{3}c_{5}c_{6}
              -\frac{5}{9}u^{6}c_{3}c_{4}c_{7}   + \cdots.    
  \end{align*}
Then,   using the second half of Lemma \ref{lem:kakikae}, the coefficients  in $(*)$ are obtained 
as follows:
\[     \lambda_{1} = 3, \; \lambda_{2} = 15, \; \lambda_{3} = 20, \;
       \lambda_{4} = \frac{10}{3}, \; \lambda_{5} = 4, \; 
      \lambda_{6} = 10, \; \lambda_{7} = 1.   
\]
Thus we have obtained 
 \begin{align*}
   I_{20}  &  \equiv  2^{27} \cdot 3^3 \cdot 5^2 \cdot 11 \cdot 17 \cdot 41 
              \left (3u^{20} +15u^{14}v + 20u^{8}v^{2} + \dfrac{10}{3}u^{2}v^{3} +4u^{10}w \right.   \\
          &  +10u^{4}vw + w^{2}  \left ) \dfrac{}{}    \right.   \\
          & \equiv  2^{27} \cdot 3^2 \cdot 5^2 \cdot 11 \cdot 17 \cdot 41
            (9u^{20}+45u^{14}v + 60u^{8}v^{2} + 10u^{2}v^{3} +12u^{10}w \\
          & +30u^{4}vw + 3w^{2}) \\
          & \equiv  2^{27} \cdot 3^2 \cdot 5^2 \cdot 11 \cdot 17 \cdot 41 
             {\tilde{I}}_{20}   \mod   \frak{b}_{20}.  
\end{align*}

Putting  
      \begin{align*}
    (**) \quad      I_{24} & \equiv 2^{32} \cdot 3^{3} \cdot 5^{2} \cdot 7 \cdot 11  \cdot 19 \cdot 199 
                    \left  ( \mu_{1}u^{24} + \mu_{2}u^{18}v + \mu_{3}u^{12}v^{2} + \mu_{4}u^{6}v^{3} \right. \\
                    & + \mu_{5}v^{4} + \mu_{6}u^{14}w + \mu_{7}u^{8}vw  + \mu_{8}u^{2}v^{2}w + \mu_{9}u^{4}w^{2} ) 
                    \mod \frak{b}_{20}
     \end{align*}
for some $\mu_{i} \in \Q$, we will proceed quite similarly. By  Lemma \ref{lem:inv.forms.E_8}  (vii), 
we have
   \begin{align*}
       I_{24} & \equiv 2^{32} \cdot 3^{3} \cdot 5^{2} \cdot 7 \cdot 11 
              \cdot 19 \cdot 199 \left ( \frac{31}{8640}c_{3}^{5}c_{4}c_{5} 
              + \frac{1}{480}c_{3}^{4}c_{5}c_{7} + \frac{337}{25920}c_{3}^{3}c_{5}^{3} \right .  \\
              &  \left. - \frac{71}{4320}c_{3}^{3}c_{4}c_{5}c_{6} +\frac{31}{240}c_{3}^{2}c_{5}c_{6}c_{7} 
              + \frac{31}{480}c_{3}c_{5}^{3}c_{6} -\frac{22}{135}c_{3}c_{4}c_{5}^{2}c_{7}    
              -\frac{1}{120}c_{4}c_{5}^{4} \right ).   
       \end{align*}
  On the other hand, in $A/(t,c_{8},\frak{a}_{20})$, we have,  for example,   
   \begin{align*}
    v^{4} & \equiv  \frac{31}{8640}c_{3}^{5}c_{4}c_{5} 
            + \frac{1}{480}c_{3}^{4}c_{5}c_{7} + \frac{337}{25920}c_{3}^{3}c_{5}^{3} 
            - \frac{71}{4320}c_{3}^{3}c_{4}c_{5}c_{6}  
            +\frac{31}{240}c_{3}^{3}c_{5}c_{6}c_{7} \\
            &+ \frac{31}{480}c_{3}c_{5}^{3}c_{6} -\frac{22}{135}c_{3}c_{4}c_{5}^{2}c_{7}    
            -\frac{1}{120}c_{4}c_{5}^{4} + \frac{11}{16}u^{7}c_{3}^{4}c_{5} - u^{7}c_{3}^{2}c_{5}c_{6}   \\
          &  -\frac{11}{18}u^{7}c_{3}c_{4}c_{5}^{2}
            + \frac{9}{160}u^{6}c_{3}^{4}c_{6} -\frac{619}{480}u^{6}c_{3}^{3}c_{4}c_{5} +  \cdots. 
    \end{align*}
Then using the second half of  Lemma \ref{lem:kakikae},  the coefficients in $(**)$
are obtained as follows:
  \begin{align*}
      & \mu_{1} = \frac{11}{5}, \; \mu_{2} = 12, \; 
        \mu_{3} = 21, \; \mu_{4} = 12, \; \mu_{5} = 1, \;
        \mu_{6} = \frac{21}{5}, \; \mu_{7} = 12, \\
       & \mu_{8} = 6, \; \mu_{9} = \frac{9}{5}.  
  \end{align*}
Thus  we have obtained 
\begin{align*}
       I_{24} & \equiv 2^{32} \cdot 3^{3} \cdot 5^{2} \cdot 7 \cdot 11 
        \cdot 19 \cdot 199  \left ( \frac{11}{5}u^{24} + 12 u^{18}v 
       + 21 u^{12}v^{2} + 12 u^{6}v^{3}  \right. \\ 
       & \left.   + v^{4} + \frac{21}{5} u^{14}w 
        + 12 u^{8}vw  + 6 u^{2}v^{2}w  + \frac{9}{5}u^{4}w^{2} \right ) \\
             & \equiv  2^{32} \cdot 3^{3} \cdot 5 \cdot 7 \cdot 11 
        \cdot 19 \cdot 199 (  11 u^{24} + 60 u^{18}v + 105 u^{12}v^{2}
                  + 60  u^{6}v^{3}  \\ 
       &   + 5 v^{4} + 21 u^{14}w 
        + 60 u^{8}vw + 30  u^{2}v^{2}w  + 9 u^{4}w^{2}) \\
       &  \equiv 2^{32} \cdot 3^{3} \cdot 5 \cdot 7 \cdot 11 
        \cdot 19 \cdot 199 \tilde{I}_{24}  \mod \frak{b}_{20}. 
 \end{align*}        

Finally,  we can  also put 
 \begin{align*}
  (***) \quad I_{30} & \equiv 2^{38} \cdot 3^{4} \cdot 5^{5} \cdot 7 \cdot 11
         \cdot 13 \cdot 61 ( \nu_{1} u^{30} + \nu_{2} u^{24}v + 
      \nu_{3} u^{18}v^{2} + \nu_{4} u^{12}v^{3}  \\
     &    + \nu_{5} u^{6}v^{4} + \nu_{6} v^{5} + \nu_{7} u^{20}w 
          + \nu_{8} u^{14}vw + \nu_{9}u^{8}v^{2}w + \nu_{10} u^{2}v^{3}w \\
     &    + \nu_{11}u^{10}w^{2} + \nu_{12}u^{4}vw^{2} + \nu_{13} w^{3})
             \mod \frak{b}_{20}      
\end{align*}
for some $\nu_{i} \in \Q$.  Then,  by Lemma \ref{lem:inv.forms.E_8}  (viii),   we have
  \begin{align*}
     I_{30}  & \equiv 2^{38} \cdot 3^{4} \cdot 5^{5} \cdot 7 \cdot 11 \cdot 13
      \cdot 61 \left ( -\frac{599}{51840}c_{3}^{5}c_{4}c_{5}c_{6} 
       + \frac{47}{34560}c_{3}^{5}c_{5}^{3} 
       + \frac{1519}{25920}c_{3}^{4}c_{5}c_{6}c_{7}  \right.   \\
    &   + \frac{6293}{7290}c_{3}^{3}c_{4}c_{5}^{2}c_{7} 
       - \frac{32537}{25920}c_{3}^{3}c_{5}^{3}c_{6}  
     + \frac{189919}{466560}c_{3}^{2}c_{4}c_{5}^{4}
      + \frac{2012}{1215}c_{3}c_{4}c_{5}^{2}c_{6}c_{7}
      - \frac{16693}{25920}c_{3}c_{5}^{4}c_{7}  \\
  &   \left.    -\frac{223}{6480}c_{4}c_{5}^{4}c_{6}
      -\frac{1}{1728}c_{5}^{6} \right ). 
  \end{align*}  
On the other hand, in $A/(t,c_{8},  \frak{a}_{20})$, we have, for example, 
   \begin{align*}
      v^{5} & \equiv \frac{31}{11520}c_{3}^{5}c_{4}c_{5}c_{6}
    -  \frac{47}{69120}c_{3}^{5}c_{5}^{3} 
     + \frac{1}{640}c_{3}^{4}c_{5}c_{6}c_{7}
     +\frac{1993}{5760}c_{3}^{3}c_{5}^{3}c_{6}
    - \frac{91}{360}c_{3}^{3}c_{4}c_{5}^{2}c_{7} \\
 &   -\frac{1279}{11520}c_{3}^{2}c_{4}c_{5}^{4}  
      -\frac{49}{135}c_{3}c_{4}c_{5}^{2}c_{6}c_{7}
    + \frac{7}{1440}c_{4}c_{5}^{4}c_{6}
    + \frac{299}{1920}c_{3}c_{5}^{4}c_{7} + \cdots,  \\
      w^{3} & \equiv  \frac{1}{162}c_{3}^{5}c_{4}c_{5}c_{6}
              -\frac{5}{81}c_{3}^{4}c_{5}c_{6}c_{7}
              +\frac{365}{648}c_{3}^{3}c_{5}^{3}c_{6}
              -\frac{1043}{2916}c_{3}^{3}c_{4}c_{5}^{2}c_{7}
              -\frac{1079}{5832}c_{3}^{2}c_{4}c_{5}^{4}  \\
            &   -\frac{226}{243}c_{3}c_{4}c_{5}^{2}c_{6}c_{7} 
                 +\frac{2}{81}c_{4}c_{5}^{4}c_{6}
              +\frac{431}{1296}c_{3}c_{5}^{4}c_{7}
              +\frac{1}{1728}c_{5}^{6} + \cdots.  
  \end{align*}
Then using the second half of Lemma \ref{lem:kakikae},  the coefficients in $(***)$
are obtained as follows:
     \begin{align*}
      &   \nu_{1} = -\frac{9}{8}, \; \nu_{2} = -3, \; \nu_{3} = 0, \; 
          \nu_{4} = -5, \; \nu_{5} = -\frac{35}{2}, \; \nu_{6} = -2, \;
          \nu_{7} = -\frac{3}{2}  \\
     &    \nu_{8} = \frac{9}{2}, \; \nu_{9} = 15, \; \nu_{10} = -5, \;
          \nu_{11} = -\frac{3}{2}, \; \nu_{12} = 3, \; \nu_{13} = -1.   
     \end{align*}
Thus we have obtained  
 \begin{align*}
     I_{30} & \equiv 2^{38} \cdot 3^{4} \cdot 5^{2} \cdot 7 \cdot 11 \cdot
    13 \cdot 61  \left ( -\frac{9}{8}u^{30}-3u^{24}v-5u^{12}v^{3}
         - \frac{35}{2}u^{6}v^{4} -2v^{5} \right.         \\
     &   \left.  -\frac{3}{2}u^{20}w  +\frac{9}{2}u^{14}vw + 15u^{8}v^{2}w 
     -5u^{2}v^{3}w  -\frac{3}{2}u^{10}w^{2} +3u^{4}vw^{2} -w^{3} \right )  \\
       & \equiv 2^{35} \cdot 3^{4} \cdot 5^{2} \cdot 7 \cdot 11 \cdot
    13 \cdot 61 ( -9u^{30}-24u^{24}v-40u^{12}v^{3}-140u^{6}v^{4}
     -16v^{5}         \\
     &    -12u^{20}w  +36u^{14}vw + 120u^{8}v^{2}w 
     -40u^{2}v^{3}w  -12u^{10}w^{2} +24u^{4}vw^{2} -8w^{3})  \\
      & \equiv 2^{35} \cdot 3^{4} \cdot 5^{2} \cdot 7 \cdot 11 \cdot
    13 \cdot 61 {\tilde{I}}_{30} \mod \frak{b}_{20}. 
     \end{align*}

Consequently,  we have established  Lemma \ref{lem:invariants}.
\end{proof}

\section{ Rational cohomology  ring of $E_{8}/T^1 \! \cdot \! E_7$}  \label{rational} 
With the above results,   we will  compute the rational cohomology  ring of $E_{8}/C_{8}$. 
We begin by  recalling  the classical results of  Borel \cite{Bor53};  
 Let $G$ be a compact connected Lie  group, $H$  a closed connected  subgroup of $G$ of maximal rank 
 and $T$ a commnon maximal torus. 
Consider  the fibration
   \[  G/H \overset{\iota}{\longrightarrow} BH  \overset{\rho}  {\longrightarrow} BG.   \] 
Since $H^*(BG;\Q)$ is a polynomial ring generated by elements of even degrees and  $H^{*}(G/H;\Q)$ 
has vanishing odd dimensional part (Hirsch  formula  \cite{Bor53}),  the Serre spectral 
sequence with rational coefficients for this fibration  collapses. 
In particular,  we have the following description of the rational cohomology  ring of $G/H$:
  \begin{align*}
       H^{*}(G/H;\Q)    \overset{\iota^{*}}{\overset{\sim}
                 {\longleftarrow}} & \, H^{*}(BH;\Q)/(\rho^{*}H^{+}(BG;\Q))      \\
                 \cong             & \, H^{*}(BT;\Q)^{W(H)}/(H^{+}(BT;\Q)^{W(G)}),  
  \end{align*}
where $H^{+} = \oplus_{i > 0} H^{i}$ and $(H^{+}(BT;\Q)^{W(G)})$ means the  ideal  of $H^*(BT;\Q)^{W(H)}$ 
generated by $H^{+}(BT;\Q)^{W(G)}$. 

We apply this result to the fibration:
   \[    E_{8}/C_{8} \overset{\iota}{\longrightarrow} BC_{8} \overset{\rho}{\longrightarrow}  BE_{8}. \]  
 Then,  using  Lemmas \ref{lem:W(C_8)-invariants}, \ref{lem:W(E_8)-invariants}, \ref{lem:invariants} and
 (\ref{eqn:generators}),  we have
     \begin{align*}
       H^{*}(E_{8}/C_{8};\Q) 
                   \cong &  \; H^{*}(BT;\Q)^{W(C_{8})}
                               /(H^{+}(BT;\Q)^{W(E_{8})})  \\ 
     \cong & \;\Q[u,J_{2},J_{6},J_{8},J_{10},J_{12},J_{14},J_{18}]
       /(I_{2},I_{8},I_{12},I_{14},I_{18},I_{20},I_{24},I_{30})  \\
      \cong & \; \Q[u,I_{2},v,I_{8},w,I_{12},I_{14},I_{18}]
       /(I_{2},I_{8},I_{12},I_{14},I_{18},I_{20},I_{24},I_{30})  \\
      \cong & \; \Q[u,v,w]/(\tilde{I}_{20},\tilde{I}_{24},
       \tilde{I}_{30}). 
     \end{align*}
Thus we have obtained the following:   
\begin{lem} \label{lem:rat.coh.E_8.C_8}
  The rational cohomology ring of $E_{8}/T^1 \! \cdot \! E_7$ is  given as follows$:$
      \[    H^{*}(E_{8}/T^1 \! \cdot \! E_7;\Q) = \Q[u,v,w] /(\tilde{I}_{20},\tilde{I}_{24}, \tilde{I}_{30}),  \]
  where  $\deg u = 2, \; \deg v = 12, \; \deg w = 20$,   
 $\tilde{I}_{20},\tilde{I}_{24}$ and $\tilde{I}_{30}$ are given by 
 $(\ref{eqn:relations})$.           
\end{lem}

\section{Integral cohomology ring of $E_{8}/T^1 \! \cdot \! E_{7}$}
\subsection{Integral cohomology ring of $E_8/T$ in low degrees}
Consider the fibration
  \[     E_{7}/T' \cong C_{8}/T \overset{i}{\longrightarrow} E_{8}/T
         \overset{p}{\longrightarrow}E_{8}/C_{8}.    \]
 Since $H^{*}(E_{8}/C_{8};\Z)$ and $H^{*}(E_{7}/T';\Z)$ have no torsion and 
vanishing odd dimensional part by Bott \cite{Bott56}, the Serre spectral sequence 
with the integral coefficients for the above fibration collapses  and the following  sequence 
  \[ \Z \rightarrow H^{*}(E_{8}/C_{8};\Z) \overset{p^{*}}
   {\rightarrow} H^{*}(E_{8}/T;\Z) \overset{i^{*}}{\rightarrow}
   H^{*}(C_{8}/T;\Z) \cong H^{*}(E_{7}/T';\Z) \rightarrow \Z  \]  
  is co-exact\footnote{This terminology is taken from \cite[\S 4]{Bau68}}, that is,  
 \begin{align*}
     &  p^{*} \text{ is injective}, \;  i^{*} \text{ is surjective  and} \\   
     &  \Ker i^{*} = ( p^{*}H^{+}(E_{8}/C_{8};\Z) ), \;  
        \text{the ideal  generated by }p^{*}H^{+}(E_{8}/C_{8};\Z).    
  \end{align*} 
Therefore we will obtain some information about the generators of 
$H^{*}(E_{8}/C_{8};\Z)$ by considering $\Ker i^{*}$. In order to investigate
$\Ker i^{*}$, we will determine $H^{*}(E_{8}/T;\Z)$ up to degrees $\leq 36$. 
First we need a simple lemma concerning the action of the cohomology operations: 
\begin{lem} \label{lem:operation}
  For the elements $t$ and $c_{i} = \sigma_{i}(t_{1},\ldots, t_{8})$ in $H^{*}(BT;\Z)$, 
we have
  \begin{align*}
   \mathrm{(i)}
     \quad & Sq^{2}(c_{2}) \equiv  c_{3} + tc_{2}, \\ 
                      & Sq^{4}(c_{3}) \equiv c_{5} + tc_{4} + c_{2}c_{3},    \\
           &  Sq^{8}(c_{5}+ tc_{4}) \equiv tc_{8} + c_{2}c_{7} + c_{3}c_{6} + c_{4}c_{5} 
              + tc_{4}^{2} + t^{2}c_{7} + t^{3}c_{6}  + t^{2}c_{2}c_{5}   + t^{2}c_{3}c_{4},   \\
           & Sq^{14}(c_{8} + c_{4}^2 + t^2 c_{6} + t^4 c_{4} + t^8) 
            \equiv  (c_{8} + t^2 c_{6} + t^4 c_{4} + t^8)(c_{7} + tc_{6}) \mod 2. \\
    \mathrm{(ii)}
          \quad & \mathcal{P}^{1}(c_{2} + 2t^2 ) \equiv c_{4} + c_{2}^{2} + t^4,   \\
                & \mathcal{P}^{3}(c_{4} + 2t^4 ) \equiv c_{5}^{2} +2 c_{4}c_{6}  +2c_{3}c_{7}
                        + 2c_{2}c_{8} + c_{3}^{2}c_{4} + c_{2}c_{4}^{2} + c_{2}^{2}c_{6}+ 2c_{2}c_{3}c_{5} \\
                        & \hspace{2.2cm}  + 2t^{10}   \mod 3.   \\
    \mathrm{(iii)} 
        \quad  &\mathcal{P}^{1}(c_{2} + t^2) \equiv c_{6} + 2c_{3}^{2} + 4c_{2}c_{4} + 2c_{2}^{3} 
                        + 2tc_{5} + tc_{2}c_{3} + 4t^{2}c_{4} + 4t^{2}c_{2}^{2} + 3t^{3}c_{3}  \\
                        & \hspace{2.1cm}  + t^{4}c_{2}  + 2t^6 \mod 5.    
  \end{align*}
\end{lem}
 
\begin{proof}
(i) follows immediately from the Wu formula: 
  \[  Sq^{2i-2}(c_{i}) \equiv   \displaystyle{\sum_{j = 0}^{i-1}}c_{2i-1-j}c_{j},  \]  
and $c_{1} = 3t \equiv t \mod 2$. 

(ii) Put $p_{i} = t_{1}^{i} +  \cdots + t_{8}^{i} \; (i \geq 0)$. 
Then $p_{i}$'s and $c_{i}$'s are related to each other by the Newton formula:
  \[     p_{n} = \sum_{i=1}^{n-1}(-1)^{i-1}p_{n-i}c_{i} + (-1)^{n-1}nc_{n}.   \]
In particular,  considering with mod 3 coefficients,  we have
   \begin{align*}
         p_{1}   &\equiv c_{1} \equiv 0, \\ 
         p_{2}   &\equiv c_{2}, \\ 
         p_{4}   &\equiv 2c_{4}+2c_{2}^{2},    \\
         p_{10}  & \equiv 2c_{5}^{2} + c_{4}c_{6} + c_{3}c_{7} + c_{2}c_{8}
                  +2c_{3}^{2}c_{4} +2c_{2}c_{4}^{2} +2c_{2}^{2}c_{6} +c_{2}^{3}c_{4}+ c_{2}c_{3}c_{5} +c_{2}^{5}. 
  \end{align*}
On the other hand,  we have 
  \begin{align*}
    \mathcal{P}^{1}(p_{2}) & \equiv \mathcal{P}^{1} \left (\sum_{i}t_{i}^{2} \right) 
      \equiv \sum_{i}\mathcal{P}^{1}(t_{i}^{2})  \equiv \sum_{i}2t_{i}
       \mathcal{P}^{1}(t_{i}) \equiv \sum_{i}2t_{i} \cdot t_{i}^{3}  \\ 
    &\equiv \sum_{i}2t_{i}^{4}   \equiv 2p_{4},    \\ 
   \mathcal{P}^{3}(p_{4}) &\equiv \mathcal{P}^{3} \left (\sum_{i}t_{i}^{4} 
   \right ) \equiv  \sum_{i}\mathcal{P}^{3}(t_{i}^{4})  \equiv 
   \sum_{i}(2 \mathcal{P}^{3}(t_{i}^{2})t_{i}^{2} +2 \mathcal{P}^{2}
   (t_{i}^{2})\mathcal{P}^{1}(t_{i}^{2}))\\
   & \equiv \sum_{i}(2 t_{i}^{6} \cdot 2t_{i}^{4}) \equiv \sum_{i}t_{i}^{10}
    \equiv p_{10}.   
 \end{align*}
Using these facts,  we have easily the required results.

(iii) Similar computation yields the required results.
\end{proof}

\begin{lem} \label{lem:int.coh.E_8/T}
       The integral cohomology ring of $E_8/T$ for degrees $\leq 36$ is given as follows$:$
    \begin{align*}
          H^{*}(E_{8}/T;\Z) & =  \Z[t_{1},\ldots,t_{8},t,\gamma_{3},\gamma_{4}, \gamma_{5},\hat{\gamma}_{6},
              \gamma_{9},  \gamma_{10},  \hat{\gamma}_{15}]  \\
     &  \hspace{0.3cm} /( \rho_{1},\rho_{2},\rho_{3},\rho_{4},\rho_{5},\hat{\rho}_{6},\rho_{8},
         \rho_{9},\rho_{10}, \rho_{12}, \rho_{14}, \hat{\rho}_{15}, \rho_{18} ),   
    \end{align*}
 where $t_{1},\ldots ,t_{8},t \in H^{2}$ are as in \S $2$, $\gamma_{i} \in H^{2i} \;
  (i = 3,4,5,6,9,10, 15)$  and
  \begin{align*}
     \rho_{1}       &= c_{1}-3t, \\
     \rho_{2}       &= c_{2}-4t^{2}, \\ 
     \rho_{3}       &= c_{3}-2\gamma_{3}, \\    
     \rho_{4}       &= c_{4}+2t^{4}-3\gamma_{4}, \\
     \rho_{5}       &= c_{5}-3t\gamma_{4} + 2t^{2}\gamma_{3} - 2\gamma_{5}, \\
     \hat{\rho}_{6} &= c_{6}-2 \gamma_{3}^{2}-t\gamma_{5}+t^{2}\gamma_{4}-2t^{6}-5\hat{\gamma}_{6},   \\
     \rho_{8}       &= -3c_{8}+3 \gamma_{4}^{2}-2\gamma_{3}\gamma_{5}+ t(2c_{7}-6\gamma_{3}\gamma_{4})
                     + t^2(2 \gamma_{3}^{2}-5\hat{\gamma}_{6}) +3t^{3}\gamma_{5} +4t^{4}\gamma_{4}-6t^{5}\gamma_{3}+t^{8},  
 \end{align*} 
 \begin{align*} 
     \rho_{9} &= 2c_{6}\gamma_{3}+tc_{8}+t^{2}c_{7}-3t^{3}c_{6}-2\gamma_{9}, \\
     \rho_{10} &=  \gamma_{5}^{2}-2c_{7}\gamma_{3}  -t^{2}c_{8}  +3t^{3}c_{7}  -3\gamma_{10}, \\
     \rho_{12} &= 15   \hat{\gamma}_{6}^2 + 2 \gamma_{3} \gamma_{4} \gamma_{5} -2 c_{7} \gamma_{5} + 2 \gamma_{3}^{4} 
                 + 10  \gamma_{3}^2 \hat{\gamma}_{6} - 3 c_{8} \gamma_{4} - 2  \gamma_{4}^3   \\
               &  + t (c_{8} \gamma_{3} -2 \gamma_{3}^2 \gamma_{5} + 4 c_{7} \gamma_{4} + 6  \gamma_{3} \gamma_{4}^{2})  
               + t^{2} (3 \gamma_{10} - 25 \gamma_{4} \hat{\gamma}_{6} - c_{7} \gamma_{3} -16  \gamma_{3}^2 \gamma_{4}) \\
               & + t^{3} (25 \gamma_{3} \hat{\gamma}_{6} - 3 \gamma_{4} \gamma_{5} + 10 \gamma_{3}^3) 
                 + t^{4} (3c_{8} + 3 \gamma_{3} \gamma_{5} + 5 \gamma_{4}^{2}) 
                 + t^{5} (-3 c_{7} - 5 \gamma_{3} \gamma_{4}) \\
               & + 4 t^{6}  \gamma_{3}^{2} - 7 t^{8} \gamma_{4} + 4 t^{9} \gamma_{3}, \\
     \rho_{14} &= c_{7}^2 - 3 c_{8} \hat{\gamma}_{6} + 6 \gamma_4 \gamma_{10} - 4c_8 \gamma_{3}^2 + 6 c_{7} \gamma_3 \gamma_4  
                  - 6  \gamma_{3}^2  \gamma_{4}^2 - 12  \gamma_{4}^2 \hat{\gamma}_6 - 2 \gamma_3 \gamma_5 \hat{\gamma}_6 \\
               &  + t(24 \gamma_3 \gamma_4 \hat{\gamma}_6 - 8c_{7} \gamma_{3}^2  - 8 c_7 \hat{\gamma}_6 + 4 c_8 \gamma_5 
                  - 6 \gamma_3 \gamma_{10} + 12 \gamma_{3}^3 \gamma_4)  \\
               &  + t^2 (-2 \gamma_3 \gamma_4 \gamma_5 + 6 \gamma_{4}^3 + 2 \gamma_{3}^2 \hat{\gamma}_6 + 20 \hat{\gamma}_{6}^2 
                          - 4 \gamma_{3}^4 - c_7 \gamma_5)  \\
               &  + t^3 (-12 \gamma_3 \gamma_{4}^2 + 8c_8 \gamma_3 - 5 c_7 \gamma_4 + 3\gamma_5 \hat{\gamma}_6) 
                  + t^4 (3 \gamma_{10} - 26 \gamma_4 \hat{\gamma}_6 + 6c_7 \gamma_3 - 4 \gamma_{3}^2 \gamma_4) \\
               &  + t^5 (24 \gamma_3 \hat{\gamma}_6 + 3 \gamma_4 \gamma_5 + 12 \gamma_{3}^3) 
                  + t^6 (-6 c_8 + 2 \gamma_{4}^2) - 4 t^7 c_7 + t^8 (6 \hat{\gamma}_6 - 6 \gamma_{3}^2) \\
               &  - 6 t^{10} \gamma_4 + 12 t^{11} \gamma_3 - 2t^{14}, \\
     \hat{\rho}_{15} &= (c_8 - t^2 c_6 + 2 t^3 \gamma_5 + 3 t^4 \gamma_4 - t^8)(c_7 - 3tc_6)  
                  - 2( \gamma_{3}^2 + c_6)(\gamma_9 - c_6 \gamma_3)  - 2\hat{\gamma}_{15}, \\
     \rho_{18}   &=  \gamma_{9}^2 - 9c_{8}\gamma_{10} - 6\gamma_{4}^2\gamma_{10} - 4\gamma_{3}^3\gamma_{9} 
                    - 10  \gamma_{3} \hat{\gamma}_{6}\gamma_{9} + 2\gamma_{3}\gamma_{5}\gamma_{10}
                    - 2\gamma_{3}\gamma_{4}\gamma_{5}\hat{\gamma}_{6} \\
               &- 6c_{7}\gamma_{3}\gamma_{4}^2 + 3c_{8}\gamma_{4}\hat{\gamma}_{6} + c_{8}\gamma_{3}^2\gamma_{4}
                   + 6\gamma_{3}^2\gamma_{4}^3 + 12\gamma_{4}^3 \hat{\gamma}_{6}+ 2c_{7}^2\gamma_{4} 
                   + 2c_{7}\gamma_{3}^2\gamma_{5} \\
              &  - 2\gamma_{3}^3\gamma_{4}\gamma_{5}+ 2c_{7}\gamma_{5}\hat{\gamma}_{6} + 4\gamma_{3}^6 
                 - 10\hat{\gamma}_{6}^3  + 18\gamma_{3}^4 \hat{\gamma}_{6} + 15 \gamma_{3}^2 \hat{\gamma}_{6}^2
                 - 9c_{7}c_{8}\gamma_{3} \\
              & + t(-2\gamma_{3}\gamma_{5}\gamma_{9}  - 24c_{7}\gamma_{4}\hat{\gamma}_{6} + 8c_{8}\gamma_{4}\gamma_{5}
                + 4c_{7}\gamma_{3}^2\gamma_{4}  + 4c_{7}\gamma_{10} - c_{8}\gamma_{9} + 2c_{7}^2 \gamma_{3}  \\
              & + 4c_{8}\gamma_{3} \hat{\gamma}_{6} + 12 \gamma_{3} \gamma_{4} \gamma_{10} 
                - 36\gamma_{3}\gamma_{4}^2 \hat{\gamma}_{6} + 12 \gamma_{3}^2\gamma_{5}\hat{\gamma}_{6} 
                + c_{8}\gamma_{3}^3 + 6\gamma_{3}^4\gamma_{5} - 18\gamma_{3}^3 \gamma_{4}^2)  \\
             &  + t^2 (24\gamma_{3}^4\gamma_{4} - 2c_{8}^2 - c_{7}\gamma_{9} - 11\gamma_{3}^2\gamma_{10} 
                   + 2\gamma_{3}\gamma_{4}\gamma_{9} - 2c_{8}\gamma_{3}\gamma_{5} + 16c_{7}\gamma_{3}\hat{\gamma}_{6} 
                   - 3c_{7}\gamma_{4}\gamma_{5}  \\
               &  + 75\gamma_{4} \hat{\gamma}_{6}^2 - 6\gamma_{4}^4 - 9c_{8}\gamma_{4}^2
                   + 81 \gamma_{3}^2\gamma_{4} \hat{\gamma}_{6} - 13 \hat{\gamma}_{6}\gamma_{10} 
                   + 4\gamma_{3}\gamma_{4}^2 \gamma_{5}  - c_{7}\gamma_{3}^3)  \\
               & + t^3(-3\gamma_{5}\gamma_{10} - 150\gamma_{3} \hat{\gamma}_{6}^2 - 135\gamma_{3}^3 \hat{\gamma}_{6}
                 + 6\gamma_{3}^2 \gamma_{9}  - 2 c_{7}\gamma_{3}\gamma_{5} + 21c_{7}\gamma_{4}^2  + 15c_{7}c_{8}  \\
                 & + 3\gamma_{4}\gamma_{5} \hat{\gamma}_{6} 
                      - 3\gamma_{3}^2 \gamma_{4} \gamma_{5} + 18\gamma_{3} \gamma_{4}^3 + 15\hat{\gamma}_{6}\gamma_{9} 
                       + 14c_{8} \gamma_{3} \gamma_{4} - 30\gamma_{3}^5)   \\
               & + t^4(-13c_{8} \hat{\gamma}_{6} + 2\gamma_{4}\gamma_{10} - 5c_{7}^2 - 33\gamma_{3}^2 \gamma_{4}^2 
                 + 3\gamma_{5} \gamma_{9} -28 \gamma_{3}\gamma_{5} \hat{\gamma}_{6} - 45\gamma_{4}^2 \hat{\gamma}_{6}  \\
               &   - 41c_{7}\gamma_{3}\gamma_{4} -13 \gamma_{3}^3 \gamma_{5} - 9c_{8} \gamma_{3}^2) \\
               & + t^5(3c_{7} \hat{\gamma}_{6} - 6\gamma_{4}^2\gamma_{5} + 23c_{7}\gamma_{3}^2  
                 + 105\gamma_{3}\gamma_{4} \hat{\gamma}_{6} - 6c_{8} \gamma_{5} - 3\gamma_{4}\gamma_{9} 
                 + 45 \gamma_{3}^3 \gamma_{4})   \\
               & + t^6(11\gamma_{4}^3 - 4\gamma_{3}\gamma_{9} + 4c_{7}\gamma_{5} + 9\gamma_{3}\gamma_{4}\gamma_{5} 
                 + 12\gamma_{3}^4 + 66\gamma_{3}^2 \hat{\gamma}_{6} + 75\hat{\gamma}_{6}^2 + 2c_{8}\gamma_{4})  \\
               & + t^7(-33\gamma_{3}\gamma_{4}^2 + 12 \gamma_{3}^2 \gamma_{5}  + 15\gamma_{5} \hat{\gamma}_{6})   
                + t^8(-4\gamma_{10} + 21\gamma_{3}^2 \gamma_{4} - 5c_{7}\gamma_{3} -3\gamma_{4} \hat{\gamma}_{6}) \\
               & + t^9(6\gamma_{9}  - 42\gamma_{3}^3 - 99\gamma_{3} \hat{\gamma}_{6}) 
                + t^{10} (-4c_{8} - 6\gamma_{4}^2 -13 \gamma_{3}\gamma_{5})  + t^{11}(3c_{7} + 27\gamma_{3}\gamma_{4})  \\
               & +t^{12} (60 \hat{\gamma}_{6} + 18\gamma_{3}^2) + 6t^{13} \gamma_{5} - 9t^{14} \gamma_{4} 
                 - 12t^{15} \gamma_{3} +10t^{18}. 
     \end{align*}
\end{lem}

\begin{proof}
   According to Toda \cite[Proposition 3.2]{Toda75}, one can give the general description of 
 $H^{*}(E_{8}/T;\Z)$  as follows:
  \begin{align*}
     H^{*}(E_{8}/T;\Z)  =  & \; \Z[t_{1},\ldots,t_{8},t,\gamma_{3},
                \gamma_{4},\gamma_{5},\gamma_{6},\gamma_{9},\gamma_{10},  \gamma_{15}] \\
              &  /(\rho_{1},\rho_{2},\rho_{3},\rho_{4},\rho_{5},\rho_{6},\rho_{8},\rho_{9},
        \rho_{10},\rho_{12},\rho_{14},\rho_{15},\rho_{18},\rho_{20},  \rho_{24},\rho_{30}), 
   \end{align*}
where $t_{1},\ldots,t_{8},t \in H^{2}$ are as above and 

 \begin{align*}
  \rho_{1} &= c_{1}-3t, \\
  \rho_{i} &= \delta_{i} - 2\gamma_{i} \;(i= 3,5,9,15),  \\
  \rho_{i} &= \delta_{i}-3\gamma_{i}\; (i = 4,10), \\ 
  \rho_{6} &= \delta_{6}-5\gamma_{6}.    
  \end{align*}
Here $\delta_{i} \; (i = 3,4,5,6,9,10,15)$ is an arbitrary element of $H^{*}(E_{8}/T;\Z)$ 
satisfying 
  \begin{align*}
 \delta_{3} &\equiv Sq^{2}(\rho_{2}) ,\; \delta_{5} \equiv Sq^{4}(\delta_{3}),
  \; \delta_{9} \equiv Sq^{8}(\delta_{5}),
  \; \delta_{15} \equiv Sq^{14}(\rho_{8})  \mod 2,  \\
 \delta_{4} &\equiv  \mathcal{P}^{1}(\rho_{2}),\; \delta_{10} \equiv 
   \mathcal{P}^{3}(\delta_{4})   \mod 3,  \\
 \delta_{6} &\equiv \mathcal{P}^{1}(\rho_{2}) \mod 5.  
  \end{align*}
Other relation $\rho_{j} \; (j = 2,8,12,14,18,20,24,30)$ is  determined by the 
maximality of the integer $n_{j}$ in 
 \begin{equation} \label{eqn:n_j}
    n_{j}  \cdot \rho_{j} \equiv \iota_{0}^{*}(I_{j})  \mod (\rho_{i} ; i < j ),          
 \end{equation}
where $\iota_{0}: E_{8}/T \longrightarrow BT$.

Now let us determine the generators and the relations explicitly; 
\begin{enumerate} 
\item In view of Lemma \ref{lem:inv.forms.E_8} (i) and  (\ref{eqn:n_j}),  we can take 
\[  \rho_{2} = c_{2} - 4t^{2}.    
\]

\item  By Lemma \ref{lem:operation} (i),  we have
\[    \delta_{3}  \equiv Sq^{2}(\rho_{2}) \equiv Sq^{2}(c_{2}) \equiv 
     c_{3} + tc_{2} \equiv c_{3} \mod (2, \rho_{1}, \rho_{2}) 
\]
and we can take $\delta_{3} = c_{3}$ so that 
\[ \rho_{3} = c_{3} - 2\gamma_{3}.   
\]

\item By Lemma \ref{lem:operation} (ii),    we have 
\[  \delta_{4}  \equiv \mathcal{P}^{1}(\rho_{2}) \equiv \mathcal{P}^{1}
  (c_{2} - 4t^{2}) \equiv c_{4} + c_{2}^{2} + t^{4}  \equiv c_{4} + 2t^{4} 
        \mod (3, \rho_{1}, \rho_{2}),    
\]
and we can take $\delta_{4} = c_{4} + 2t^{4} $ so that 
\[ \rho_{4} = c_{4} + 2t^{4} - 3\gamma_{4}.    
\]

\item By Lemma \ref{lem:operation} (i),  we have 
\begin{align*}
\delta_{5} & \equiv Sq^{4}(\delta_{3}) \equiv Sq^{4}(c_{3}) \equiv 
     c_{5} + tc_{4} + c_{2}c_{3}  \equiv c_{5} + tc_{4}  \mod (2, \rho_{1}, \rho_{2},\rho_{3}, \rho_{4})  \\
    & \equiv c_{5} - 3t\gamma_{4} + 2t^{2}\gamma_{3}   
                           \mod (2, \rho_{1}, \rho_{2}, \rho_{3}, \rho_{4}),  
\end{align*}
and we can take $\delta_{5} = c_{5} - 3t\gamma_{4} + 2t^{2}\gamma_{3} $ 
so that 
\[  \rho_{5} = c_{5}-3t\gamma_{4}+2t^{2}\gamma_{3}-2\gamma_{5}.   
\]

\item By Lemma \ref{lem:operation} (iii),   we have 
   \begin{align*}
      \delta_{6} &\equiv  \mathcal{P}^{1}(\rho_{2}) \equiv \mathcal{P}^{1}
  (c_{2} - 4t^{2}) \equiv  c_{6} + 2c_{3}^{2} + 4c_{2}c_{4}
     + 2c_{2}^{3} + 2tc_{5} + tc_{2}c_{3} \\
     &+ 4t^{2}c_{4} + 4t^{2}c_{2}^{2} + 3t^{3}c_{3} + t^{4}c_{2}  + 2t^6 \mod (5, \rho_{1})  \\
     & \equiv  c_{6} + 3\gamma_{3}^2 + 4t\gamma_{5} + t^2 \gamma_{4}+ 3t^6 \mod (5, \rho_{1}, 
       \rho_{2}, \rho_{3}, \rho_{4}, \rho_{5})  \\ 
     &  \equiv  c_{6} - 2\gamma_{3}^{2} - t\gamma_{5} + t^{2}\gamma_{4} -2t^{6}    \mod(5, \rho_{1}, 
      \rho_{2}, \rho_{3}, \rho_{4}, \rho_{5}),   
   \end{align*}
and we can take $\delta_{6} = c_{6} - 2\gamma_{3}^{2} - t\gamma_{5}  + t^{2}\gamma_{4}  -2t^{6}$ so that 
\[ \hat{\rho}_{6} = c_{6} - 2\gamma_{3}^{2} - t\gamma_{5} + t^{2}\gamma_{4} -2t^{6} -5\hat{\gamma}_{6}.    
\]

\item  By Lemma \ref{lem:inv.forms.E_8}  (ii), we have 
\begin{equation}  \label{eqn:I_8} 
 \begin{array}{ll}   
  I_{8} &\equiv 2^{14} \cdot 3 \cdot 5 
        \{ -18c_{8}-3c_{3}c_{5}+2c_{4}^{2}
            +  t(12c_{7}-3c_{3}c_{4}) + t^2(-6c_{6}+3c_{3}^{2})   \medskip  \\
        &+12t^{3}c_{5}+2t^{4}c_{4}-12t^{5}c_{3}+14t^{8}\} \mod (I_2).  \medskip 
\end{array}
\end{equation}   
On the other hand, by the relations $\rho_{3}$, $\rho_{4}$, $\rho_{5}$, $\hat{\rho}_{6}$, 
we have 
\begin{equation}  \label{eqn:c_{3}..c_{6}}
\begin{array}{rl} 
     c_{3} &= 2\gamma_{3}, \medskip   \\
     c_{4} &= 3\gamma_{4}-2t^{4},  \medskip \\
     c_{5} &= 2\gamma_{5}+3t\gamma_{4}-2t^{2}\gamma_{3},  \medskip   \\
     c_{6} &= 5\hat{\gamma}_{6} + 2 \gamma_{3}^{2} + t\gamma_{5}-t^{2}\gamma_{4} + 2t^{6}.    
\end{array}      
\end{equation} 
Substituting (\ref{eqn:c_{3}..c_{6}}) into (\ref{eqn:I_8}), we have  
 \begin{align*} 
     I_{8}   & \equiv 2^{15} \cdot 3^{2} \cdot 5 
          \{ -3c_{8}+3\gamma_{4}^{2}  -2\gamma_{3}\gamma_{5} + t(2c_{7} - 6\gamma_{3}\gamma_{4})
          +  t^2 (2\gamma_{3}^{2} - 5\hat{\gamma}_{6}) \\
        & +3t^{3}\gamma_{5} + 4t^{4}\gamma_{4}
          -6t^{5}\gamma_{3}+t^{8} \}    \mod (\rho_{1}, \rho_{2}, \rho_{3}, \rho_{4}, \rho_{5}, \hat{\rho}_{6}).    
  \end{align*}
Hence, by (\ref{eqn:n_j}), we have   
\[   2^{15} \cdot 3^2 \cdot 5 \, \rho_{8} \equiv \iota_{0}^{*}(I_8) 
          \mod (\rho_{1}, \rho_{2}, \rho_{3}, \rho_{4}, \rho_{5}, \hat{\rho}_{6}) 
\] 
and it follows the form of $\rho_{8}$.

\item  By Lemma \ref{lem:operation} (i) ,  we have 
 \begin{align*}
 \delta_{9} & \equiv  Sq^{8}(\delta_{5}) \equiv Sq^{8}(c_{5} + tc_{4})
    \equiv tc_{8} + c_{2}c_{7} + c_{3}c_{6} + c_{4}c_{5} + tc_{4}^{2} + t^2c_{7}  \\
           &  + t^3c_{6} + t^2c_{2}c_{5} + t^2c_{3}c_{4}  \mod (2, \rho_{1}) \\
           & \equiv tc_{8} + t^2c_{7} + t^3c_{6}  
              \mod (2,\rho_{1}, \rho_{2}, \rho_{3}, \rho_{5})  \\
           & \equiv tc_{8} + t^{2}c_{7} -3t^{3}c_{6} +2c_{6}\gamma_{3} 
              \mod (2,\rho_{1}, \rho_{2}, \rho_{3}, \rho_{5}), 
 \end{align*}   
and we can take $\delta_{9} = 2c_{6}\gamma_{3} +  tc_{8} + t^{2}c_{7} -3t^{3}c_{6}$ so that 
\[ \rho_{9} =  2c_{6}\gamma_{3} +  tc_{8} + t^{2}c_{7} -3t^{3}c_{6}   -2\gamma_{9}.    
\]

\item By Lemma \ref{lem:operation} (ii),  we have 
    \begin{align*}
        \delta_{10} &\equiv \mathcal{P}^{3}(\delta_{4}) \equiv \mathcal{P}^{3}   (c_{4}+2t^{4}) 
                     \equiv c_{5}^{2} + 2c_{4}c_{6} + 2c_{3}c_{7} + 2c_{2}c_{8} + c_{3}^{2}c_{4}    \\  
                    &+ c_{2}c_{4}^{2} + c_{2}^{2}c_{6} +  2c_{2}c_{3}c_{5} + 2t^{10}  \mod(3, \rho_{1})     \\
                    & \equiv \gamma_{5}^2 + c_{7}\gamma_{3} + 2t^2 c_{8} 
                      \mod (3,\rho_{1}, \rho_{2}, \rho_{3}, \rho_{4}, \rho_{5}, \hat{\rho}_{6}) \\
                    & \equiv {\gamma_{5}}^{2} - 2c_{7}\gamma_{3}  -t^{2}c_{8}  + 3t^{3}c_{7}
                      \mod (3, \rho_{1}, \rho_{2}, \rho_{3}, \rho_{4}, \rho_{5}, \hat{\rho}_{6}), 
          \end{align*}
and we can take $\delta_{10} = \gamma_{5}^{2} - 2c_{7}\gamma_{3} -t^{2}c_{8} + 3t^{3}c_{7}$ so that 
\[   \rho_{10} = \gamma_{5}^{2}-2c_{7}\gamma_{3} -t^{2}c_{8} +  3t^{3}c_{7} - 3\gamma_{10}.    
\]
    
\item  By (\ref{eqn:inv.forms.E_8}), (\ref{eqn:Newton}) and (\ref{eqn:d_n}),  we obtain  
  \begin{align*}
      I_{12} & \equiv 2^{18} \cdot 3^4 \cdot 5 \cdot 7 
                \left \{    \dfrac{3}{5}{c_6}^2 - c_{5}c_{7} - c_{4}c_{8} + \dfrac{1}{6}c_{3}c_{4}c_{5} 
                       - \dfrac{2}{27} c_{4}^3 - \dfrac{1}{10}c_{3}^2 c_{6}  + \dfrac{1}{40} c_{3}^4  \right. \\
             & + t \left (\dfrac{1}{2}c_{3}c_{8} +  \dfrac{7}{3}c_{4}c_{7} - \dfrac{3}{5}c_{5}c_{6}  
                    - \dfrac{1}{5}c_{3}^2 c_{5}  + \dfrac{1}{6} c_{3}c_{4}^2 \right )  \\
             & + t^2 \left (\dfrac{2}{5}c_{5}^2 - \dfrac{5}{2} c_{3}c_{7} - \dfrac{2}{3} c_{4}c_{6} 
                            - \dfrac{1}{6} c_{3}^2c_{4}  \right )  \\
             & + t^3 \left ( \dfrac{19}{10} c_{3}c_{6} - \dfrac{2}{3} c_{4}c_{5} - \dfrac{1}{5}c_{3}^3 \right ) 
                   + t^4 \left ( - \dfrac{1}{9}c_{4}^2  + \dfrac{19}{30} c_{3}c_{5}  \right )
                   + t^5 \left ( \dfrac{14}{3} c_{7} + \dfrac{1}{2}c_{3}c_{4} \right )  \\
             & \left.  + t^6 \left ( - \dfrac{56}{15} c_{6} + \dfrac{23}{30} c_{3}^2  \right )  
                  - \dfrac{2}{15} t^7 c_{5} - \dfrac{5}{9} t^8 c_{4} - \dfrac{22}{15}t^9 c_{3} 
                  + \dfrac{154}{135} t^{12}  \right \}    \mod (I_2) \\
             & \equiv  2^{18} \cdot 3^4 \cdot 5 \cdot 7  
                \{  15 {\hat{\gamma}_{6}}^2 + 2 \gamma_3 \gamma_4 \gamma_5 -2 c_7 \gamma_5 + 2 \gamma_{3}^4 
                   + 10 \gamma_{3}^2 \hat{\gamma}_6  - 3 c_8 \gamma_4 - 2 \gamma_{4}^3   \\
               &  + t (c_8 \gamma_3 -2 \gamma_{3}^2 \gamma_5 + 4 c_7 \gamma_4 + 6 \gamma_3 \gamma_{4}^2)  
                  + t^2 (3 \gamma_{10} - 25 \gamma_4 \hat{\gamma}_6 - c_7 \gamma_{3} -16 \gamma_{3}^2 \gamma_4) \\
               &  + t^3 (25 \gamma_3 \hat{\gamma}_6 - 3 \gamma_4 \gamma_5 + 10 \gamma_{3}^3) 
                  + t^4 (3c_8 + 3 \gamma_3 \gamma_5 + 5 \gamma_{4}^2) 
                  + t^5 (-3 c_7 - 5 \gamma_{3} \gamma_{4}) \\
               &  + 4 t^6 \gamma_{3}^2 - 7 t^8 \gamma_{4} + 4 t^9 \gamma_{3} \}  
                  \mod (\rho_{1}, \rho_{2}, \rho_{3},  \rho_{4}, \rho_{5}, \hat{\rho}_{6}, \rho_{10}).
  \end{align*}  
Hence,  we have 
\[  2^{18} \cdot 3^4 \cdot 5 \cdot 7 \; \rho_{12} \equiv \iota_{0}^{*}(I_{12}) 
    \mod (\rho_{1}, \rho_{2}, \rho_{3}, \rho_{4}, \rho_{5}, \hat{\rho}_{6}, \rho_{10}).    
\]
From this and (\ref{eqn:n_j}),  the form of $\rho_{12}$ follows.    
Quite similarly,  we have 
      \begin{align*} 
         2^{20} \cdot 3^2 \cdot 5^2 \cdot 7 \cdot 11 \; \rho_{14}  & \equiv \iota_{0}^{*}(I_{14})  
                            \mod (\rho_{1}, \rho_{2}, \rho_{3}, \rho_{4}, \rho_{5}, \hat{\rho}_{6}, \rho_{10}), \\  
         2^{26} \cdot 3^4 \cdot 5^2 \cdot 7 \cdot 13 \; \rho_{18} & \equiv \iota_{0}^{*}(I_{18}) 
                            \mod (\rho_{1}, \rho_{2}, \rho_{3}, \rho_{4}, \rho_{5}, \hat{\rho}_{6}, \rho_{9}, \rho_{10}),  
       \end{align*}
and it follows the forms of $\rho_{14}$ and $\rho_{18}$.

\item  Finally,  we will determine the relation $\hat{\rho}_{15}$. Since 
      \begin{align*} 
           \rho_{8} & = -3c_{8}+3{\gamma_{4}}^{2}-2\gamma_{3}\gamma_{5}+ t(2c_{7}-6\gamma_{3}\gamma_{4})
                    + t^2(2{\gamma_{3}}^{2}-5 \hat{\gamma}_{6}) +3t^{3}\gamma_{5} + 4t^{4}\gamma_{4}   \\
                    & -6t^{5}\gamma_{3}+t^{8} \\
                    & \equiv c_{8} + c_{4}^2 + t^2 c_{6} + t^4 c_{4} + t^8 \mod (2, \rho_{1}, \rho_{2}, \rho_{3}, 
                       \rho_{4}, \rho_{5}, \hat{\rho}_{6}), 
      \end{align*}
 we have, by Lemma \ref{lem:operation} (i), 
  \begin{align*} 
        \delta_{15} & \equiv Sq^{14}(\rho_{8})   \\
                    & \equiv Sq^{14}(c_{8} + c_{4}^2 + t^2c_{6} + t^4 c_{4} + t^8)  \\
                    & \equiv (c_{8} + t^2c_{6} + t^4c_{4} + t^8) (c_{7} + tc_{6})  
                        \mod (2, \rho_{1}, \rho_{2}, \rho_{3}, \rho_{4}, \rho_{5}, \hat{\rho}_{6}) \\
                    & \equiv (c_{8} - t^2 c_{6} + 2t^3 \gamma_{5} + 3t^4\gamma_{4} - t^8)(c_{7} - 3tc_{6}) 
                            - 2(c_{6} + \gamma_{3}^2 )(\gamma_{9} - c_{6}\gamma_{3})  \\
                    &  \mod (2, \rho_{1}, \rho_{2}, \rho_{3}, \rho_{4}, \rho_{5}, \hat{\rho}_{6}, \rho_{9})
  \end{align*}
and we can take
\[   \delta_{15} = (c_{8} - t^2 c_{6} + 2t^3 \gamma_{5} + 3t^4\gamma_{4} - t^8)(c_{7} - 3tc_{6}) 
                            - 2(c_{6} + \gamma_{3}^2 )(\gamma_{9} - c_{6}\gamma_{3})   
\] 
so that 
\[ \hat{\rho}_{15}   =  (c_{8} - t^2 c_{6} + 2t^3 \gamma_{5} + 3t^4\gamma_{4} - t^8)(c_{7} - 3tc_{6}) 
                            - 2(c_{6} + \gamma_{3}^2 )(\gamma_{9} - c_{6}\gamma_{3})  - 2 \hat{\gamma}_{15}.  
\] 
\end{enumerate} 
Consequently,  we have established the lemma.
\end{proof}

\begin{rem}
 Since $H^*(E_{8}/T;\Z)$ has a free $\Z$-basis consisting of Schubert classes, 
 our generators $\gamma_{3}, \gamma_{4}, \gamma_{5}, \hat{\gamma}_{6}, \gamma_{9}, \gamma_{10}, \hat{\gamma}_{15}$
 can be expressed as certain $\Z$-linear combinations of Schubert classes. 
 The precise expression is given in \cite{Kaji-Nak2} $($see also \cite{DZ08}$)$. 
\end{rem}

In order to determine $\Ker i^*$,  we need the  result on the integral cohomology ring $H^*(E_7/T';\Z)$, that 
was  computed by the author. The result is restated as follows:
\begin{thm}[\cite{Nak01}, Theorem 5.9]   \label{thm:E_7/T}
  The integral cohomology ring of   $E_{7}/T'$ is given as follows$:$
       \begin{align*}
            H^{*}(E_{7}/T';\Z)  & = \Z[t_{1}', \ldots ,t_{7}',t',
                   \gamma_{3}',\gamma_{4}',\gamma_{5}',\gamma_{9}']        \\
           &  \hspace{0.3cm}   /(\rho_{1}',\rho_{2}',\rho_{3}',\rho_{4}',\rho_{5}',\rho_{6}',
                    \rho_{8}',\rho_{9}',\rho_{10}', \rho_{12}', \rho_{14}', \rho_{18}'),          
        \end{align*}
where $t_{1}', \ldots, t_{7}',t' \in H^{2}$ are as in $\S 2$, $\gamma_{i}' \in H^{2i} \; (i = 3,4,5,9)$ and 
   \begin{align*}
      \rho_{1}' &= c_{1}'- 3t',\\ 
      \rho_{2}' &= c_{2}' - 4{t'}^{2}, \\
      \rho_{3}' &= c_{3}' - 2\gamma_{3}', \\ 
      \rho_{4}' &= c_{4}' + 2{t'}^{4} - 3\gamma_{4}',  \\
      \rho_{5}' &= c_{5}' -3{t'}\gamma_{4}' + 2{t'}^{2}\gamma_{3}' - 2\gamma_{5}',  \\
      \rho_{6}' &= {\gamma_{3}'}^{2} + 2c_{6}' -2t'\gamma_{5} -3{t'}^{2}\gamma_{4}'+ {t'}^{6},   \\    
      \rho_{8}' &= 3{\gamma_{4}'}^{2} - 2\gamma_{3}'\gamma_{5}' + t'(2c_{7}' -6 \gamma_{3}'\gamma_{4}')
                   -9{t'}^{2}c_{6}' 
                   +12{t'}^{3}\gamma_{5}' + 15{t'}^{4}\gamma_{4}' - 6{t'}^{5}\gamma_{3}'-{t'}^{8},   \\
      \rho_{9}' &= 2c_{6}'\gamma_{3}' +{t'}^{2}c_{7}' -3{t'}^{3}c_{6}' -2\gamma_{9}',  \\
      \rho_{10}'&= {\gamma_{5}'}^{2} -2c_{7}'\gamma_{3}' +3{t'}^{3}c_{7}', \\
      \rho_{12}' &= 3c_{6}'^2 - 2{\gamma_{4}'}^3 - 2c_{7}' \gamma_{5}' + 2\gamma_{3}'\gamma_{4}'\gamma_{5}'  
                  + t'(4c_{7}'\gamma_{4}' - 2c_{6}'\gamma_{5}' + 6\gamma_{3}'{\gamma_{4}'}^2)  \\ 
               &  + t'^2(-3c_{7}'\gamma_{3}' + 3c_{6}'\gamma_{4}') 
                  + t'^3(-12\gamma_{4}'\gamma_{5}' + 5c_{6}'\gamma_{3}') 
                  + t'^4(-2\gamma_{3}'\gamma_{5}' - 15{\gamma_{4}'}^2)  \\
                & - 10t'^6c_{6}' + 12t'^7\gamma_{5}'+ 19t'^8 \gamma_{4}' - 6t'^9 \gamma_{3}' - 2t'^{12}, \\ 
    \rho_{14}'  &= {c_{7}'}^2 + 6c_{7}'\gamma_{3}'\gamma_{4}' - 2c_{6}'\gamma_{3}'\gamma_{5}' - t'^2c_{7}'\gamma_{5}'
                  + t'^3(-9c_{7}'\gamma_{4}' + 3c_{6}'\gamma_{5}') - 6t'^4 c_{7}'\gamma_{3}'   + 9t'^7 c_{7}', \\
   \rho_{18}'   &= -{\gamma_{9}'}^2 + 2c_{6}'c_{7}'\gamma_{5}' + 6c_{7}'\gamma_{3}'{\gamma_{4}'}^2 - 2{c_{7}'}^2 \gamma_{4}' 
                   -2c_{6}'\gamma_{3}'\gamma_{4}'\gamma_{5}' + 2c_{6}'\gamma_{3}'\gamma_{9}' \\
                &  + t'(-6{c_{7}'}^2\gamma_{3}' + 24c_{6}'c_{7}'\gamma_{4}') 
                  + t'^2(-25c_{7}'\gamma_{4}'\gamma_{5}' + c_{7}'\gamma_{9}' - 18c_{6}'c_{7}'\gamma_{3}') \\
                & + t'^3(-45c_{7}'{\gamma_{4}'}^2 + 20c_{7}'\gamma_{3}'\gamma_{5}' + 3c_{6}'\gamma_{4}'\gamma_{5}' 
                   - 3c_{6}'\gamma_{9}')  \\
                & + t'^4(11c_{7}'^2 + 2c_{6}'\gamma_{3}'\gamma_{5}' + 48c_{7}'\gamma_{3}'\gamma_{4}') + 51t'^5 c_{6}'c_{7}'
                  - 53t'^6 c_{7}'\gamma_{5}'   \\ 
               &  + t'^7 (-69c_{7}'\gamma_{4}' - 3c_{6}'\gamma_{5}') + 16t'^8 c_{7}'\gamma_{3}' + 15t'^{11} c_{7}'. 
  \end{align*}
\end{thm}

\begin{rem}
In order to deal with the higher relations $\rho_{12}', \rho_{14}'$ and $\rho_{18}'$, 
the  author     make   use of the integral cohomology ring of the homogeneous space 
$E_7/T^1 \! \cdot \! \mathit{Spin}(12)$ $($\cite{Nak01}$)$. 
Here  we expressed the relations $\rho_{12}'$, $\rho_{14}'$ and $\rho_{18}'$
  in terms of the generators $t_{1}', \ldots, t_{7}', t', \gamma_{3}', \gamma_{4}', \gamma_{5}', \gamma_{9}'$
  by computing the $W(E_{7})$-invariants $I_{12}'$, $I_{14}'$ and $I_{18}'$  explicitly $($see $\ref{W(C_8)-invariants})$.  
\end{rem}

By Lemma \ref{lem:int.coh.E_8/T} and Theorem \ref{thm:E_7/T}, we can find the generators of the ideal $\Ker i^{*}$: 
\begin{prop} \label{cor:kerneli^*}
   For the induced homomorphism  
    \[  i^{*} : H^{*}(E_{8}/T;\Z) \longrightarrow H^{*}(C_{8}/T;\Z) \cong  H^{*}(E_{7}/T';\Z),   \]
we obtain that 
    \[ \Ker i^{*} = (u, \tilde{\gamma}_{6},\gamma_{10}, \hat{\gamma}_{15}),    \]
where $\tilde{\gamma}_{6} = 2 \hat{\gamma}_{6}+\gamma_{3}^{2}-t^{2}\gamma_{4}+t^{6}$.
\end{prop}
\begin{proof}
 By (\ref{eqn:t_i}),  we have 
  \begin{equation} \label{eqn:i*-images1}
     i^{*}(t_{i}) = t_{i}' \; (1 \leq i \leq 7), \; i^{*}(t_{8}) = 0, \;
     i^{*}(t) = t',   
  \end{equation}
and therefore
   \begin{equation} \label{eqn:i*-images2}
    i^{*}(c_{n}) = c_{n}' \; (1 \leq n \leq 7), \; i^{*}(c_{8}) = 0.  
   \end{equation}
Then it is verified directly that 
\begin{equation} \label{eqn:i*-images3}
  \begin{array}{cl}
      i^{*}(\gamma_{i}) & = \gamma_{i}' \; (i = 3,4,5,9), \medskip  \\ 
      i^{*}(\hat{\gamma}_{6}) &= c_{6}'-t'{\gamma_{5}}' -{t'}^{2}{\gamma_{4}}',  \medskip \\
      i^{*}(\gamma_{10})& = 0,  \medskip  \\  
      i^{*}(\hat{\gamma}_{15})&= 0. 
  \end{array}
\end{equation}

Now we put 
  \[ I = (u, \tilde{\gamma}_{6}, \gamma_{10}, \hat{\gamma}_{15}), \] 
the ideal of $H^{*}(E_{8}/T;\Z)$  generated by the elements in the parenthesis. 
Using (\ref{eqn:i*-images1}), (\ref{eqn:i*-images3}) and Theorem \ref{thm:E_7/T}, 
we see that $I$ is contained in $\Ker i^{*}$.   Hence there is an induced map 
   \[   H^{*}(E_{8}/T;\Z)/I \longrightarrow H^{*}(C_{8}/T;\Z) \cong  H^{*}(E_{7}/T';\Z).   \]
Then,  by Lemma \ref{lem:int.coh.E_8/T}, we have 
   \begin{align*}
   H^{*}(E_{8}/T;\Z)/I  \cong & \;  \Z[t_{1},\ldots, t_{8},t,\gamma_{3},
                                  \gamma_{4},\gamma_{5}, \hat{\gamma}_{6},\gamma_{9},\gamma_{10}, \hat{\gamma}_{15}] \\
                              & /(u, \rho_{1},\rho_{2},\rho_{3},\rho_{4},\rho_{5}, \hat{\rho}_{6},
                                 \tilde{\gamma}_{6},\rho_{8},\rho_{9},\rho_{10},\gamma_{10}, \rho_{12}, \rho_{14}, \hat{\rho}_{15},
                                 \hat{\gamma}_{15}, \rho_{18}) \\
                       \cong & \;\Z[t_{1},\ldots,t_{7},t,\gamma_{3},\gamma_{4},\gamma_{5}, \gamma_{9}] \\
                             & /(\rho_{1},\rho_{2},\rho_{3},\rho_{4},\rho_{5},\tilde{\gamma}_{6}, \rho_{8},\rho_{9},\rho_{10}, 
                               \rho_{12}, \rho_{14}, \rho_{18}) 
    \end{align*}
for degrees $\leq 36$. Then  it follows from  Lemma \ref{lem:int.coh.E_8/T}, Theorem \ref{thm:E_7/T}, 
(\ref{eqn:i*-images1}), (\ref{eqn:i*-images2}) and (\ref{eqn:i*-images3})  that 
      \begin{align*}
             i^*(\rho_{i}) & \equiv \rho_{i}' \; (i = 1, 2, 3, 4, 5, 8, 9, 10, 12, 14, 18), \\
             i^*(\tilde{\gamma}_{6})  & \equiv \rho_{6}'.  
      \end{align*}  
Therefore this map induces  an isomorphism  and the assertion follows.
\end{proof}

\subsection{Generators of $H^*(E_8/T^1 \! \cdot \! E_7;\Z)$}  \label{u,v,w,x}
From Proposition \ref{cor:kerneli^*},  we see that $H^{*}(E_{8}/C_{8};\Z)$ is generated  as a ring by some four
elements $\tilde{u} \in H^{2},\; \tilde{v} \in H^{12}, \; \tilde{w}\in  H^{20}$ and $\tilde{x} \in H^{30}$
 such that 
     \begin{equation} \label{eqn:ideals}
 (\tilde{u},\tilde{v},\tilde{w}, \tilde{x}) = (u, \tilde{\gamma}_{6}, \gamma_{10}, \hat{\gamma}_{15})  
   \end{equation} 
as ideals. So our next task is to  describe  these generators in the ring $H^{*}(E_8/T;\Z)$. 
 Hereafter we  identify  $H^{*}(E_{8}/C_{8};\Z)$ with the subalgebra $\Im p^{*}$ of $H^{*}(E_{8}/T;\Z)$.

Firstly, by Lemma \ref{lem:J6,J10} and (\ref{eqn:bar{c_i}}),  we have 
  \begin{equation} \label{eqn:tilde(J_6)}
  \begin{array}{cl}
   \tilde{J}_{6} &:= \dfrac{1}{2^{10} \cdot 3^{2} \cdot 5}  J_{6}   \medskip \\
                 & \equiv \dfrac{2}{5}c_{6}+ \dfrac{1}{20}c_{3}^{2}+c_{5}
                   \left  (-\dfrac{1}{5}t -\dfrac{1}{2}u \right )
                   +c_{4} \left (\dfrac{1}{2}tu + \dfrac{1}{3}u^{2} \right )  \medskip \\ 
                 &   + c_{3} \left (-\dfrac{1}{5}t^{3}     
                    - \dfrac{1}{2}t^{2}u - \dfrac{1}{2}u^{3} \right ) 
                  +  \dfrac{1}{5}t^{6} +  t^{5}u - \dfrac{1}{3}t^{4}u^{2}   
                   + t^{3}u^{3}+ t^{2}u^{4}- tu^{5}   \medskip \\ 
               & +\dfrac{273}{640}u^{6} \mod (I_{2}),    
   \end{array}
\end{equation}

 \begin{equation} \label{eqn:tilde(J_10)} 
  \begin{array}{cl}
    \tilde{J}_{10} & :=  \dfrac{1}{2^{14} \cdot 3^{3} \cdot 5 \cdot 7} J_{10}  \medskip \\
                   & \equiv    \dfrac{1}{12}c_{5}^{2} - \dfrac{1}{3}c_{3}c_{7}
                     + \dfrac{1}{2}c_{3}c_{6}u + c_{4}c_{5} \left (-\dfrac{1}{6}t 
                      - \dfrac{1}{6}u \right )
                     + c_{3}c_{5} \left (\dfrac{1}{6}t^{2}-\dfrac{23}{84}u^{2} \right )  \medskip  \\
                    & + c_{4}^{2} \left (\dfrac{1}{12}t^{2} + \dfrac{1}{6}tu + 
                     \dfrac{1}{14}u^{2} \right)    
                    + c_{3}^{2} \left (\dfrac{1}{12}t^{4} + \dfrac{23}{84}t^{2}u^{2} 
                     -\dfrac{37}{96}u^{4} \right )   \medskip  \\
                     & + c_{3}c_{4} \left (-\dfrac{1}{6}t^{3}
                   - \dfrac{1}{6}t^{2}u -\dfrac{23}{84}tu^{2} + \dfrac{1}{3}u^{3} \right )  \medskip  \\
                   & +c_{7} \left (t^{3}+\dfrac{1}{6}t^{2}u-\dfrac{17}{42}tu^{2}
                     +\dfrac{5}{14}u^{3} \right )   
                     +c_{6} \left (-\dfrac{3}{2}t^{3}u - \dfrac{3}{14}t^{2}u^{2} 
                     + \dfrac{1}{2}tu^{3}  -\dfrac{37}{84}u^{4} \right )   \medskip \\
                    & + c_{5}\left (-\dfrac{1}{3}t^{5} + \dfrac{1}{6}t^{4}u 
                     + \dfrac{23}{21}t^{3}u^{2} - \dfrac{1}{3}t^{2}u^{3} 
                     + \dfrac{1}{24}tu^{4} +\dfrac{71}{336}u^{5} \right )  \medskip \\
                    & + c_{4} \left (\dfrac{1}{3}t^{6} + \dfrac{1}{6}t^{5}u 
                      + \dfrac{17}{42}t^{4}u^{2}
                      - \dfrac{2}{3}t^{3}u^{3} - \dfrac{1}{3}t^{2}u^{4}+ \dfrac{5}{16}tu^{5}
                      -\dfrac{131}{504}u^{6} \right )  \medskip \\ 
                    & + c_{3} \left (-\dfrac{1}{3}t^{7}  
                      + \dfrac{1}{6}t^{6}u -\dfrac{23}{21}t^{5}u^{2}- \dfrac{1}{3}t^{4}u^{3}
                      + \dfrac{37}{24}t^{3}u^{4} + \dfrac{73}{48}t^{2}u^{5}
                      - \dfrac{3}{2}tu^{6} + \dfrac{239}{336}u^{7} \right )  \medskip  \\
                     & + \dfrac{1}{3}t^{10} - \dfrac{1}{3}t^{9}u + \dfrac{7}{6}t^{8}u^{2}
                       + \dfrac{2}{3}t^{7}u^{3} - \dfrac{7}{8}t^{6}u^{4} - \dfrac{35}{8}t^{5}u^{5}
                       + \dfrac{233}{72}t^{4}u^{6} + \dfrac{7}{24}t^{3}u^{7}  \medskip  \\
                     & -\dfrac{129}{56}t^{2}u^{8} + \dfrac{215}{168}tu^{9} 
                       - \dfrac{208601}{645120}u^{10} \mod (I_{2}).     \medskip   
            \end{array}
\end{equation}
On the other hand,  by Lemma \ref{lem:int.coh.E_8/T}, we have 
\begin{equation}  \label{eqn:c_3..c_6}
\begin{array}{rl} 
     c_{3} &= 2\gamma_{3}, \medskip   \\
     c_{4} &= 3\gamma_{4}-2t^{4},  \medskip \\
     c_{5} &= 2\gamma_{5}+3t\gamma_{4}-2t^{2}\gamma_{3},  \medskip   \\
     c_{6} &= 5\hat{\gamma}_{6} + 2 \gamma_{3}^{2} + t\gamma_{5}-t^{2}\gamma_{4} + 2t^{6}    \medskip  
\end{array}      
\end{equation} 
in $H^*(E_8/T;\Z)  \hookrightarrow H^{*}(E_{8}/T;\Q)$. Substituting (\ref{eqn:c_3..c_6}) 
into (\ref{eqn:tilde(J_6)}), we have 
  \begin{equation} \label{eqn:tilde(J_6).2}
  \begin{array}{cl}
   \tilde{J}_{6} & = \dfrac{1}{2^{10} \cdot 3^2 \cdot 5}J_{6}  \medskip \\
                 & = 2\hat{\gamma}_{6}+  \gamma_{3}^{2} - u\gamma_{5} +\gamma_{4}(-t^{2}+u^{2})
                    - u^{3}\gamma_{3}+ t^{6} -t^{4}u^{2}  + t^{3}u^{3}  \medskip  \\
                 &  + t^{2}u^{4}-tu^{5}+ \dfrac{273}{640}u^{6}.  \medskip 
  \end{array}
\end{equation} 
Similarly,  substituting (\ref{eqn:c_3..c_6}) into (\ref{eqn:tilde(J_10)}) 
and  using the relations $\rho_{8},\rho_{9}$ and $\rho_{10}$, we have
\begin{equation} \label{eqn:tilde(J_10).2}
  \begin{array}{cl}
  \tilde{J}_{10} &= \dfrac{1}{2^{14} \cdot 3^{3} \cdot 5 \cdot 7}J_{10}  \medskip \\
                 & =  \gamma_{10} + u\gamma_{9} - u^{3}c_{7} 
                     -u\gamma_{4}\gamma_{5}  + 2u^{2}{\gamma}_{4}^{2} - 2u^{2}\gamma_{3}\gamma_{5} 
                  +  \gamma_{3}\gamma_{4}(-6tu^{2}+2u^{3})  \medskip \\
                 & + \gamma_{3}^{2} \left (2t^{2}u^{2}+2tu^{3}+ \dfrac{7}{24}u^{4} \right)  
                   + \hat{\gamma}_{6} \left (-5t^{2}u^{2}+5tu^{3}+\dfrac{55}{12}u^{4} \right )  \medskip \\
                &  + \gamma_{5} \left (t^{4}u + 3t^{3}u^{2} + t^{2}u^{3} 
                        - \dfrac{55}{24}u^{5} \right )  \medskip \\  
                  &    + \gamma_{4} \left (6t^{4}u^{2} - 3t^{3}u^{3} - \dfrac{103}{24}t^{2}u^{4} 
                       - tu^{5} + \dfrac{79}{24}u^{6} \right )  \medskip   \\
                  &    + \gamma_{3} \left (-6t^{5}u^{2}-2t^{4}u^{3}+4t^{3}u^{4}
                  +6t^{2}u^{5}-4tu^{6}  - \dfrac{31}{24}u^{7} \right )  \medskip  \\
                & +4t^{7}u^{3}+\dfrac{55}{24}t^{6}u^{4}-6t^{5}u^{5}-\dfrac{7}{24}t^{4}u^{6}
                  + \dfrac{79}{24}t^{3}u^{7} + \dfrac{31}{24}t^{2}u^{8} 
                  - \dfrac{55}{24}tu^{9}  \medskip \\ 
                & + \dfrac{666919}{645120}u^{10}.   \medskip   
   \end{array}
\end{equation}

Now let us determine our generators $\tilde{u},\tilde{v}, \tilde{w}$ and $\tilde{x}$. 
Obviously,   we can take $u = t_{8}$ as our generator $\tilde{u}$.

Next, since $H^*(E_8/C_8;\Q)$ is generated by $u, J_{6}$ and  $J_{10}$ (see \S \ref{rational}),   we can  put 
    \[ \tilde{v} = \alpha \tilde{J}_{6} + \beta u^{6}    \]
for some $\alpha, \beta \in \Q$. On the other hand,  by (\ref{eqn:ideals}),  we can express  
    \[ \tilde{v} = \tilde{\gamma}_{6} + f = 2 \hat{\gamma}_{6} + \gamma_{3}^{2} - t^2 \gamma_{4} + t^{6}  + f  \]
for some element $f \in (u) \cap H^{12}(E_{8}/T;\Z)$. Then,  using (\ref{eqn:tilde(J_6).2}),  we have 
  \begin{align*}
    &  2\hat{\gamma}_{6}+  \gamma_{3}^{2}-t^{2}\gamma_{4}+t^{6}+ f 
       =   \alpha(2\hat{\gamma}_{6}+  \gamma_{3}^{2}-t^{2}\gamma_{4}+t^{6}) \\
    & + \alpha (-u\gamma_{5}+u^{2}\gamma_{4}-u^{3}\gamma_{3}-t^{4}u^{2}
      + t^{3}u^{3} + t^{2}u^{4} - tu^{5}) 
                      + \left (\frac{273}{640}\alpha + \beta \right )u^{6},  
 \end{align*} 
and we can  take   $  \alpha = 1, \; \beta = -\dfrac{273}{640}$ . 
Thus we see that 
  \begin{equation} \label{eqn:generator.v}
   \begin{array}{cl}
  v & = \dfrac{1}{2^{10} \cdot 3^{2} \cdot 5} J_{6} - \dfrac{273}{640}u^{6}  \medskip \\
    & = 2\hat{\gamma}_{6}+ \gamma_{3}^{2}-u\gamma_{5}+\gamma_{4}(-t^{2}+u^{2})
        -u^{3}\gamma_{3} + t^{6} - t^{4}u^{2} + t^{3}u^{3} + t^{2}u^{4} - tu^{5}   \medskip 
   \end{array}
  \end{equation}
can be chosen as our generator $\tilde{v}$.

Similarly,  we can  put
  \[  \tilde{w} = \lambda \tilde{J}_{10} + \mu u^{4}v + \nu u^{10}   \]
for some $\lambda, \mu, \nu \in \Q$. On the other hand,  by (\ref{eqn:ideals}),  we can express 
 \[ \tilde{w} = \gamma_{10} + g  \]
for some element $g \in (u, \tilde{\gamma}_{6}) \cap H^{20}(E_{8}/T;\Z)$. 
Then,  using (\ref{eqn:tilde(J_10).2}),  we can  take $\lambda = 1$ and hence
  \begin{align*}
      \gamma_{10} + g & =  \gamma_{10} + u\gamma_{9} - u^{3}c_{7}  
                          -u\gamma_{4}\gamma_{5} + 2u^{2} \gamma_{4}^{2} 
                          -2u^{2}\gamma_{3}\gamma_{5} + \gamma_{3}\gamma_{4}(-6tu^{2}+2u^{3}) \\
                      &  +  \gamma_{3}^{2}  \left \{ 2t^{2}u^{2}+2tu^{3}+ \left (\mu 
                       + \frac{7}{24} \right )u^{4} \right \}  \\
                        &  + \hat{\gamma}_{6} \left \{ -5t^{2}u^{2}+5tu^{3}
                         + \left(2\mu +  \frac{55}{12} \right )u^{4} \right \} \\
                     & +\gamma_{5} \left \{ t^{4}u + 3t^{3}u^{2} + t^{2}u^{3} + 
                         \left (-\mu - \frac{55}{24} \right )u^{5} \right \}  \\
                     & +\gamma_{4} \left  \{ 6t^{4}u^{2} - 3t^{3}u^{3} 
                       + \left (-\mu - \frac{103}{24} \right )t^{2}u^{4}- tu^{5} + 
                       \left  (\mu + \frac{79}{24} \right )u^{6} \right  \}  \\
                     &  + \gamma_{3} \left \{ -6t^{5}u^{2} - 2t^{4}u^{3} 
                        + 4t^{3}u^{4} + 6t^{2}u^{5}  -4tu^{6} +  \left (-\mu - \frac{31}{24}
                        \right )  u^{7}   \right \}   \\
                     & + 4t^{7}u^{3}+ \left (\mu + \frac{55}{24} \right )t^{6}u^{4} -6t^{5}u^{5} 
                       + \left (-\mu -\frac{7}{24} \right )t^{4}u^{6} + \left (\mu
                       + \frac{79}{24} \right )t^{3}u^{7} \\
                    & + \left (\mu + \frac{31}{24} \right )t^{2}u^{8} 
                      + \left (-\mu- \frac{55}{24}\right )tu^{9} 
                      + \left (\nu + \frac{666919}{645120} \right )u^{10},   
       \end{align*}
and we can take $ \mu = -\dfrac{55}{24}, \; \nu = -\dfrac{666919}{645120}$.
Thus we see that 
\begin{equation} \label{eqn:generator.w}
  \begin{array}{cl}
    w & = \dfrac{1}{2^{14} \cdot 3^{3} \cdot 5 \cdot 7}J_{10} 
- \dfrac{55}{24}u^{4}v - \dfrac{666919}{645120}u^{10}  \medskip  \\
    & = \gamma_{10} + u\gamma_{9} - u^{3}c_{7}  -u\gamma_{4}\gamma_{5} +2u^{2} \gamma_{4}^{2}
        -2u^{2}\gamma_{3}\gamma_{5} + \gamma_{3}\gamma_{4}( -6tu^{2} + 2u^{3})  \medskip \\
    &+ \gamma_{3}^{2}(2t^{2}u^{2} + 2tu^{3} - 2u^{4})  + \hat{\gamma}_{6}(-5t^{2}u^{2} + 5tu^{3})
     + \gamma_{5}(t^{4}u + 3t^{3}u^{2} + t^{2}u^{3})  \medskip \\
   & +\gamma_{4}(6t^{4}u^{2} - 3t^{3}u^{3} -2t^{2}u^{4} -tu^{5} + u^{6}) \medskip \\
   & +\gamma_{3}(-6t^{5}u^{2} - 2t^{4}u^{3} + 4t^{3}u^{4} + 6t^{2}u^{5} - 4tu^{6} + u^{7})  \medskip \\
   & + 4t^{7}u^{3} -6t^{5}u^{5} + 2t^{4}u^{6}  + t^{3}u^{7} - t^{2}u^{8}  \medskip 
 \end{array}
\end{equation}
can be chosen as our generator $\tilde{w}$.

Finally,   we have to find an element $x$ of degree $30$ such that 
$x \equiv \hat{\gamma}_{15} \mod (u, v, w)$ in $H^*(E_{8}/T;\Z)$. 
In order to identify the element $x$, we make use of the result on the  mod $2$ cohomology ring $H^*(E_{8}/T;\Z/2\Z)$ 
(up to degrees $\leq 30$)  
due to Kono and Ishitoya (\cite[Theorem 3.10]{Kono-Ishi87})\footnote{Note that 
our generators $\gamma_{5}, \gamma_{9}, \hat{\gamma}_{15}$ are slightly different from 
$\gamma_{5}, \gamma_{9}, \gamma_{15}'$ in \cite{Kono-Ishi87}.}. 
Using their result, it can be checked directly that $u^{15} \equiv 0$ in $H^*(E_{8}/T;\Z/2\Z)$.
This means that $u^{15}$ is divisible by $2$ in the ring $H^*(E_{8}/T;\Z)$. Therefore
there exists an element $x \in H^{30}(E_{8}/T;\Z)$ such that $2x = u^{15}$. 
Explicit form of $x$ in $H^*(E_{8}/T;\Z)$ is given in Appendix \ref{degree15}.  
Furthermore,   we can check directly that 
      \[  x \equiv \hat{\gamma}_{15} \mod (u, v, w)    \]
(see also Appendix \ref{degree15}).  Hence the element $x$ can be chosen as our generator $\tilde{x}$.

\subsection{Integral cohomology ring of $E_8/T^1 \! \cdot \! E_7$}
Using the element $x$,  we can rewrite $\tilde{I}_{30}$ as follows: 
  \begin{align*} 
        \tilde{I}_{30} & =  \;  -9u^{30} - 24u^{24}v  -12u^{20}w  + 36u^{14}vw -40u^{12}v^3  -12u^{10}w^2  + 120u^8 v^2w \\
                       &  -140u^{6}v^{4} + 24u^{4}vw^{2} -40u^2v^3w -16v^5 -8w^{3} \\
                       &=  \;  -36x^2 - 48u^{9}vx  -24u^{5}wx  + 36u^{14}vw -40u^{12}v^3  -12u^{10}w^2  + 120u^8 v^2w \\
                       &  -140u^{6}v^{4} + 24u^{4}vw^{2} -40u^2v^3w -16v^5 -8w^{3} \\
                       &=  \; 4(-9x^2 - 12u^{9}vx  -6u^{5}wx  + 9u^{14}vw -10u^{12}v^3  -3u^{10}w^2  + 30u^8 v^2w \\
                         &  -35u^{6}v^{4} + 6u^{4}vw^{2} - 10u^2v^3w -4v^5 -2w^{3}).  
  \end{align*}
Therefore, in view of Lemma \ref{lem:rat.coh.E_8.C_8},  we obtain the following  main result of 
this paper. 
\begin{thm}  \label{thm:E_8/T^1E_7}
   The integral cohomology ring of $E_{8}/T^1 \! \cdot \! E_{7}$ is given as follows$:$
          \[   H^*(E_{8}/T^1 \! \cdot \! E_{7};\Z) = \Z [u, v, w, x]/(r_{15}, r_{20}, r_{24}, r_{30}),  \] 
 where $\deg u = 2$, $\deg v = 12$, $\deg w = 20$, $\deg x = 30$ and  
  \begin{align*}
     r_{15}   &= u^{15} - 2x, \\
     r_{20}   &= 9u^{20} + 45u^{14}v + 12u^{10}w + 60u^8v^2 + 30u^4vw + 10u^2v^3 + 3w^2, \\
     r_{24}   &=  11u^{24} + 60u^{18}v +  21u^{14}w + 105u^{12}v^2 + 60u^8vw + 60u^{6}v^3 + 9u^4 w^2 \\
              & + 30 u^2 v^2 w + 5v^4,   \\
     r_{30}   &= -9x^2 - 12u^{9}vx  -6u^{5}wx  + 9u^{14}vw -10u^{12}v^3  -3u^{10}w^2  + 30u^8 v^2w \\
              &  -35u^{6}v^{4} + 6u^{4}vw^{2} - 10u^2v^3w -4v^5 -2w^{3}.
\end{align*}   
\end{thm}

\begin{rem} 
As remarked in the introduction, the integral cohomology ring $H^*(E_{8}/T^1 \! \cdot \! E_{7};\Z)$ is also 
computed by Duan and Zhao  in terms of Schubert classes $($\cite[Theorem 7]{DZ05}$)$. 
 We will verify that both results  completely coincide in our forthcoming paper \cite{Kaji-Nak1}.  
\end{rem}

\section{Integral cohomology of $E_8/E_7$} 
In order to determine the integral cohomology  of $E_8/E_7$, we consider the Gysin exact sequence 
associated with the  following circle bundle 
   \begin{equation} \label{eqn:circle_bundle} 
       S^1  \longrightarrow E_8/E_7 \overset{\pi}{\longrightarrow} E_{8}/T^1 \! \cdot \! E_7,   
    \end{equation}    
  where $\pi$ is the natural projection. In this case, it reduces to the following short exact 
  sequeces:
   \begin{equation} \label{eqn:Gysin}
        \begin{array}{lll}
          0 \longrightarrow H^{\rm{odd}}(E_8/E_7;\Z)  & \longrightarrow H^*(E_8/T^1 \! \cdot \! E_7;\Z) \medskip \\ 
          & \overset{\cup u}{\longrightarrow} H^{*}(E_8/T^1 \! \cdot \! E_7;\Z) 
            \overset{\pi^*}{\longrightarrow} H^{\rm{even}}(E_8/E_7;\Z) \longrightarrow 0,   
  \end{array}
\end{equation}  
where $H^{\rm{even}} = \oplus_{i \geq 0} H^{2i}$ and $H^{\rm{odd}} = \oplus_{i \geq 0} H^{2i+1}$. 
From the exactness of (\ref{eqn:Gysin}),  it follows that $H^{\rm{even}}(E_8/E_7;\Z)$ is 
isomorphic to $H^{*}(E_8/T^1 \! \cdot  E_7;\Z)/(u)$. Define the elements
 $z_{i} \; (i = 12, 20, 30)$ of $H^*(E_8/E_7;\Z)$ as 
    \[  z_{12} = \pi^{*}(v), \; z_{20} = \pi^{*}(w), \; z_{30} = \pi^{*}(x).    \] 
Then,   by Theorem \ref{thm:E_8/T^1E_7},  we obtain 
  \begin{align*}
        H^{\text{even}} (E_{8}/E_{7};\Z)  &  =  \Z[z_{12}, z_{20}, z_{30}]/
                                                  (2z_{30}, 3 z_{20}^2, 5 z_{12}^4, 
                                                  4 z_{12}^5 + 2 z_{20}^3 + 9 z_{30}^2)     \\
                                           & =  \Z[z_{12}, z_{20}, z_{30}]/
                                                  (2z_{30}, 3 z_{20}^2, 5 z_{12}^4, 
                                                  z_{12}^5 + z_{20}^3 + z_{30}^2).  
  \end{align*} 
By Poincar\'{e} duality,  there exist elements $z_{i} \in H^{i}(E_8/E_7;\Z) \; (i = 59$, $71$, $79$, $83$, 
$91$, $95$, $103$, $115)$ such that 
   \[  z_{12}^3 z_{20} z_{59} = z_{12}^2 z_{20} z_{71} = z_{12}^3 z_{79} = z_{12}z_{20}z_{83} 
        = z_{12}^2 z_{91} = z_{20}z_{95} = z_{12}z_{103} = z_{115}.  \]        
 Then it is not hard to show that 
    \begin{align*} 
      z_{71} &=  z_{12}z_{59}, \\
      z_{79} &= z_{20}z_{59}, \\  
      z_{83} &= z_{12}^2 z_{59}, \\
      z_{91} &= z_{12}z_{20}z_{59}, \\
      z_{103} &= z_{12}^2 z_{20} z_{59}, \\
      z_{115} &= z_{12}^3 z_{20} z_{59}. 
     \end{align*}              

Summing up the results,  we obtain the following:           
  \begin{cor}   \label{cor:E_8/E_7} 
    The structure of $H^{*}(E_{8}/E_{7};\Z)$ is given by the following table$:$
     \begin{center}  
      \begin{tabular}{|l|c|c|}
      \noalign{\hrule height0.8pt}
      \hfil     $H^{k}(E_{8}/E_{7};\Z)$ &  {\rm  elements}    \\
      \hline
               $H^{0}  \cong \Z$  &  $1$     \\  
      \hline
               $H^{12} \cong \Z$ &  $z_{12}$  \\  
      \hline
               $H^{20} \cong \Z$ &  $z_{20}$   \\    
      \hline 
               $H^{24} \cong \Z$ & $z_{12}^2$  \\
      \hline
               $H^{30} \cong \Z/2\Z$ & $z_{30}$   \\
      \hline
               $H^{32} \cong \Z$    & $z_{12}z_{20}$ \\ 
      \hline
               $H^{36} \cong \Z$   & $z_{12}^3$ \\
      \hline
               $H^{40} \cong \Z/3\Z$ & $z_{20}^2$    \\
      \hline
               $H^{42} \cong \Z/2\Z$ & $z_{12}z_{30}$   \\
      \hline  
               $H^{44} \cong \Z$     & $z_{12}^2 z_{20}$ \\
      \hline  
               $H^{48} \cong \Z/5\Z$  & $z_{12}^4$  \\
     \hline 
               $H^{50} \cong \Z/2\Z$  & $z_{20}z_{30}$  \\
     \hline  
               $H^{52} \cong  \Z/3\Z$  & $z_{12} z_{20}^2$  \\
     \hline   
               $H^{54} \cong \Z/2\Z$  & $z_{12}^2 z_{30}$  \\
     \hline
               $H^{56} \cong \Z$      & $z_{12}^3 z_{20}$   \\   
     \hline  
               $H^{59} \cong \Z$      & $z_{59}$  \\
      \hline 
              $H^{62} \cong \Z/2\Z$  & $z_{12}z_{20}z_{30}$  \\
      \hline
              $H^{64} \cong \Z/3\Z$  & $z_{12}^2 z_{20}^2$ \\
      \hline 
              $H^{66} \cong \Z/2$  & $z_{12}^3 z_{30}$  \\
      \hline  
              $H^{68} \cong \Z/5\Z$  & $z_{12}^4 z_{20}$  \\
      \hline
               $H^{71} \cong \Z$      & $z_{12}z_{59}$  \\
      \hline
               $H^{74} \cong \Z/2\Z$  & $z_{12}^2 z_{20} z_{30}$ \\
      \hline
               $H^{76} \cong \Z/3\Z$  & $z_{12}^3 z_{20}^2$ \\
      \hline  
               $H^{79} \cong \Z$      & $z_{20}z_{59}$   \\
      \hline 
               $H^{83} \cong \Z$      & $z_{12}^2 z_{59}$ \\
      \hline  
             $H^{86} \cong \Z/2\Z$   & $z_{12}^3 z_{20} z_{30}$ \\
      \hline 
               $H^{91} \cong \Z$      & $z_{12}z_{20}z_{59}$  \\
      \hline 
               $H^{95} \cong \Z$      & $z_{12}^3 z_{59}$ \\
      \hline  
               $H^{103} \cong \Z$     & $z_{12}^2 z_{20} z_{59}$  \\
      \hline
               $H^{115} \cong \Z$     & $z_{12}^3 z_{20} z_{59}$   \\       
   \noalign{\hrule height0.8pt}
   \end{tabular} 
   \end{center} 
\end{cor}

\section{Appendix}
\subsection{Generator of degree $15$}   \label{degree15}
In \ref{u,v,w,x},  we defined the element $x$ as a rational cohomology class given by 
  \[  x = \dfrac{1}{2} u^{15}  \quad \text{in} \quad H^{30}(E_8/T;\Q).   \]
We need to show that $x$ is in fact an  integral cohomology class.  By Lemma \ref{lem:int.coh.E_8/T}, 
the following relations hold in $H^{*}(E_8/T\; \Z)$: 
\begin{equation}  \label{eqn:relations1}
\begin{array}{cll}
  c_{1} \! &= 3t, \medskip \\
  c_{2} \! &= 4t^2,\medskip \\
  c_{3} \! &= 2\gamma_{3}, \medskip  \\
  c_{4} \! &= 3\gamma_{4} - 2t^4, \medskip \\
  c_{5} \! &= 2\gamma_{5} + 3t\gamma_{4} - 2t^2 \gamma_{3}, \medskip \\
  c_{6} \! &= 5\hat{\gamma}_{6} + 2\gamma_{3}^2 + t\gamma_{5} - t^2 \gamma_{4} + 2t^6. \medskip 
\end{array}
\end{equation}  
  Note that, by   (\ref{eqn:c_{8}}) and (\ref{eqn:relations1}), the following relation holds:
  \begin{equation} \label{eqn:c_8} 
  \begin{array}{cl}  
     c_{8} & = uc_{7} - 5u^2 \hat{\gamma}_{6} - 2u^2 \gamma_{3}^2 + \gamma_{5}(-tu^2 + 2u^3)  
              + \gamma_{4} (t^2 u^2 + 3tu^3 - 3u^4)  \medskip \\
           &  + \gamma_{3} (-2t^2 u^3 + 2u^5)  
            -2t^6 u^2 + 2t^4 u^4 - 4t^2 u^6 + 3tu^7 - u^8.   \medskip 
  \end{array} 
  \end{equation}   
  Using (\ref{eqn:c_8}), we can rewrite the higher relations 
 $\rho_{8}, \rho_{9}, \rho_{10}, \rho_{12}, \rho_{14}$ and $\hat{\rho}_{15}$. For example, 
  \begin{equation*} 
   \begin{array}{cll} 
       \rho_{8} &= -3c_{8} + 3\gamma_{4}^2 - 2\gamma_{3}\gamma_{5} + t(2c_{7} - 6\gamma_{3}\gamma_{4}) 
                 + t^2(2\gamma_{3}^2 - 5\hat{\gamma}_{6}) + 3t^3 \gamma_{5}  + 4t^4 \gamma_{4}  \medskip \\
                &- 6t^5 \gamma_{3} + t^8   \medskip   \\
                &= 3\gamma_{4}^2 -2\gamma_{3}\gamma_{5} +  c_{7} (2t - 3u) -6t\gamma_{3}\gamma_{4} 
                 + \hat{\gamma}_{6} (-5t^2 + 15u^6)    \medskip \\
                &+ \gamma_{3}^2 (2t^2 + 6u^2 )  +  \gamma_{5} (3t^3 + 3tu^2 - 6u^3)  
                 + \gamma_{4} (4t^4 -3t^2 u^2 -9tu^3 + 9u^4)    \medskip \\
                &+ \gamma_{3} (-6t^5 + 6t^2 u^3 - 6u^5)  + t^8 + 6t^6 u^2 -6t^4 u^4 + 12 t^2 u^6 - 9tu^7 + 3u^8. 
   \end{array}
 \end{equation*}  
Using the relations $\rho_{i} \; (i = 1, 2, 3, 4, 5, 8, 9, 10, 12, 14)$, $\hat{\rho}_{6}$,  $\hat{\rho}_{15}$,  we can  rewrite 
the element $x = \dfrac{1}{2}u^{15}$ as follows:
  \begin{align*}
     x &= \dfrac{1}{2}u^{15}  \\
       &= \dfrac{1}{2} \left \{u^{15} - \hat{\rho}_{15} + u\rho_{14} - u^3 \rho_{12}  +  (3t^4u - 3t^2 u^3)\rho_{10}  \right. \\
                       & \left. + (\hat{\gamma}_{6} + (t + u) \gamma_{5} + u^2 \gamma_{4} + u^3 \gamma_{3}  
                                + t^6 + t^4 u^2 + t^3 u^3 + t^2 u^4 + tu^5 + u^6)\rho_{9}  - 39u^7 \rho_{8} \right \}  \\
       & + (-tu \gamma_{3} + 2t^3 u^2 - 5tu^4) \rho_{10}  
 \end{align*}
 \begin{align*} 
  &=  \hat{\gamma}_{15} - 20 \gamma_{3} \hat{\gamma}_{6}^2 + 3 \gamma_{3}^{2}\gamma_{9}  
     -23 \gamma_{3}^3 \hat{\gamma}_{6} 
   -6 \gamma_{3}^{5}  + 4 \hat{\gamma}_{6}\gamma_{9}  \\
  &+ 3 u \gamma_{4}\gamma_{10} - u \gamma_{5}\gamma_{9}  
    -3 u \gamma_{3}^{2} \gamma_{4}^{2} +3 uc_{7}\gamma_{3}\gamma_{4}  -6 u \gamma_{4}^2 \hat{\gamma}_{6} 
    + (-3t + 2 u ) \gamma_{3}^{3} \gamma_{5}     +(-4 t  + 4u) \gamma_{3}\gamma_{5} \hat{\gamma}_{6}  \\
  & +(-t^2 - u^2) \gamma_{4}\gamma_{9}  +(t^2 + tu -u^2 )c_{7} \gamma_{3}^{2}  
  +(9t^2+ 12 tu + 5u^2) \gamma_{3}\gamma_{4} \hat{\gamma}_{6}  \\
  & +  (5t^2 +6tu + 2 u^2) \gamma_{3}^{3} \gamma_{4}  + (3 t^2 +4tu + u^2)  c_{7} \hat{\gamma}_{6}        \\
  & +(-6t^3 -2t^2u - 6tu^2 + 5u^3) \gamma_3^{4}    - u^3 \gamma_{3}\gamma_{9} 
    + (3t^2 u + u^3 ) \gamma_{4}^{3}  + (2 t^2u + 3tu^2) c_{7} \gamma_{5}  \\
  &   +(-45t^3 + 10t^2 u -40 tu^2) \hat{\gamma}_{6}^{2} + (t^3 -2t^2 u + tu^2 - u^3 )\gamma_{3} \gamma_{4} \gamma_{5}  \\
  & +(-33 t^3 + t^2u - 31tu^2 + 13u^3) \gamma_{3}^{2} \hat{\gamma}_{6} \\ 
  & +( -2 t^4 - 4 t^3u -3t u^3 + 3u^4) c_{7} \gamma_{4}    
    +(-9 t^4 -6t^3u -18t^2u^2 + 5tu^3 -3u^4)\gamma_{5} \hat{\gamma}_{6}  \\  
  & +(-3t^4 -3t^3 u -7t^2 u^2 + 5tu^3 -4u^4) \gamma_{3}^{2}\gamma_{5} 
    +(-t^4 -6 t^3u  -t^2u^2 -3tu^3) \gamma_{3} \gamma_{4}^2    \\
   & + (-3t^4 u -6 t^3 u^2 + 3t^2 u^3  +15 t u^4  )\gamma_{10}   
     + (-3t^4 u +  t^3 u^2 + 5 t^2 u^3 + 10tu^4  -u^5 )c_{7} \gamma_{3}   \\
  &  +(15 t^5 -2t^4u + 3t^3 u^2 + 14t^2 u^3 -16tu^4 + 3u^5) \gamma_{3}^2 \gamma_{4}  \\
  & + (39t^5 - 13t^4u + 8t^3 u^2 + 35 t^2 u^3 -31tu^4 - 3u^5)\gamma_{4} \hat{\gamma}_{6}     \\
  & +(t^6 -t^4u^2 - t^3 u^3 - t^2 u^4 - tu^5 -u^6)\gamma_{9}    \\
  & +(- 13t^6 + 12t^5 u +5t^4 u^2 -56t^3 u^3 + 8t^2 u^4 +21tu^5 + 2u^6) \gamma_{3} \hat{\gamma}_{6}       \\
  & + ( 6 t^6 + 3t^5 u +2t^4 u^2 +  7 t^3u^3  + t^2 u^4 -8tu^5 + 3u^6)\gamma_{4}  \gamma_{5}    \\
  & +(-8 t^6 + 6t^5u + 2t^4u^2 -22t^3 u^3 + 6t^2u^4 +8tu^5 -2u^6) \gamma_{3}^{3}   \\
  & + (-6t^7 + t^6 u -7t^4 u^3  + 5 t^3 u^4  + 3t^2 u^5 + 3tu^6 - 63 u^7) \gamma_{4}^{2}   \\
  & + (-t^7 +2t^6 u + t^5 u^2 -11 t^4 u^3 + 6t^3 u^4 + 5t^2 u^5 +6tu^6 + 39u^7)\gamma_{3}\gamma_{5}   \\
  & + (2t^8 + 6t^7 u +3t^6u^2 -4t^5 u^3  - 15 t^4 u^4  + 6t^3 u^5 + 3t^2 u^6  -40tu^7 + 59 u^8)c_{7}   \\
  & + (3t^8 +  t^6 u^2 +11t^5 u^3 + 14t^4 u^4 -20t^3 t^5 -4 t^2 u^6 +118tu^7 + 3u^8 )\gamma_{3} \gamma_{4}    \\  
  & + (- 48 t^9 + 3 t^8u  -41 t^7 u^2  + 18 t^6 u^3 + 16 t^5 u^4 -13t^4 u^5  - 67t^3 u^6 + 125t^2 u^7  \\
  & - 15tu^8  -291u^9) \hat{\gamma}_{6}     \\ 
  & +(-18t^9 -3t^8 u -16t^7 u^2 + 10t^6u^3 -4t^5u^4 -8t^4 u^5 -16t^3 u^6  -23t^2 u^7 -10tu^8  -115u^9) \gamma_{3}^{2}  \\
  & +(-6t^{10} -3t^9 u -9t^8 u^2 + 5t^7 u^3 -5 t^6 u^4  -14t^4 u^6 -52t^3 u^7 + 6t^2 u^8 
    -60tu^9  +117 u^{10})\gamma_{5}   \\
  & + (18t^{11} -3t^{10} u+ 5t^9 u^2 + 11t^8 u^3 - 28t^7 u^4 + 8t^6 u^5 + 20t^5 u^6  -64t^4 u^7 
    -15t^3 u^8 \\
  &  + 54t^2 u^9 + 178t u^{10} - 177u^{11}) \gamma_{4} \\
  & +(-2t^{12} +6t^{11}u + 2t^{10}u^2 -20 t^9 u^3 +11t^8 u^4 + 22t^7 u^5 -8t^6 u^6 + 83t^5 u^7  \\
  & + 15t^4 u^8 + 5t^3 u^9  -116 t^2 u^{10}  
    + tu^{11} + 117 u^{12}) \gamma_{3}        \\
 &  -12t^{15} - t^{14}u -10t^{13}u^2 + 6 t^{12} u^3 + 7t^{11}u^4 -13t^{10}u^5  -31t^9 u^6  +9 t^8 u^7 -t^7 u^8 \\
 &  -118t^6 u^9 -18t^5 u^{10} + 131t^4 u^{11} -6t^3 u^{12}   - 233t^2u^{13}  
    + 175tu^{14} -58 u^{15},       
\end{align*} 
which has shown that $x$ is in fact an integral cohomology class. 

Next,  we have to  show that 
      \[ x \mod (u, v, w ) \equiv \hat{\gamma}_{15}.  \]
By (\ref{eqn:generator.v}) and (\ref{eqn:generator.w}),  we have 
  \begin{align*}
     v & \equiv 2 \hat{\gamma}_{6} + \gamma_{3}^2 - t^2 \gamma_{4} + t^6 \mod (u), \\
     w & \equiv \gamma_{10} \mod (u, v).   
  \end{align*} 
Therefore, in the ring $H^*(E_8/T;\Z)/(u, v, w)$, the following  relations hold:
\begin{equation} \label{eqn:mod(u,v,w)}
  \begin{array}{cl}
     u & = 0,  \medskip \\
     \gamma_{3}^2 & = -2 \hat{\gamma}_6 + t^2 \gamma_{4} - t^6, \medskip \\
     \gamma_{10} & = 0. \medskip
  \end{array}
\end{equation}   
On the other hand, we determined  the ring $H^*(E_8/T;\Z)$ up to degrees $\leq 36$ (Lemma \ref{lem:int.coh.E_8/T}). 
Taking (\ref{eqn:mod(u,v,w)}) into account, we can show directly that $x = \hat{\gamma}_{15}$ in the 
ring $H^*(E_8/T;\Z)/(u, v, w)$.

\subsection{Relations of $H^*(E_8/T;\Z)$ in degrees $\leq 18$}
  In Lemma \ref{lem:int.coh.E_8/T}, we took  the element $\hat{\gamma}_{6}$ as one of the 
  ring generators of $H^*(E_{8}/T;\Z)$, so that  the relation 
 of degree $12$ is given as follows:
   \begin{equation*} 
     \hat{\rho}_{6} = c_{6} - 2\gamma_{3}^2 - t\gamma_{5} + t^2 \gamma_{4} - 2t^6 - 5\hat{\gamma}_{6}. 
   \end{equation*}   
Since $\tilde{\gamma}_{6} = 2 \hat{\gamma}_{6} + \gamma_{3}^2 - t^2 \gamma_{4} + t^6$  (see Proposition
\ref{cor:kerneli^*}),  we have $\hat{\gamma}_{6} = -2 \tilde{\gamma}_{6} + c_{6} - t\gamma_{5} - t^2 \gamma_{4}$, 
and  we can replace $\hat{\gamma}_{6}$ with $\tilde{\gamma}_{6}$ as a new generator, 
so that the relation of degree $12$ changes to 
 \begin{equation*} 
    \tilde{\rho}_{6}  = \gamma_{3}^2 + 2c_{6} - 2t\gamma_{5} - 3t^2 \gamma_{4} + t^6 - 5\tilde{\gamma}_{6}.  
 \end{equation*} 

Using the element $\tilde{\gamma}_{6}$, other relations $\rho_{8}, \rho_{12}, \rho_{14}, \rho_{18}$ 
change to the following:
\begin{align*} 
  \rho_{8}  &= 3\gamma_{4}^2 -2 \gamma_{3}\gamma_{5} -3 c_{8} + t(2c_{7} - 6\gamma_{3}\gamma_{4}) 
             + t^2 (-9c_{6} + 20\tilde{\gamma}_{6})  + 12t^3 \gamma_{5} + 15t^4 \gamma_{4} \\
            & - 6t^5 \gamma_{3} - t^8, \\
  \rho_{12} &= 3c_{6}^2 - 2\gamma_{4}^3 -2c_{7}\gamma_{5} + 2\gamma_{3}\gamma_{4}\gamma_{5}  
             + 10\tilde{\gamma}_{6}^2-10c_{6}\tilde{\gamma}_{6} - 3c_{8}\gamma_{4}  \\
            &+ t(4c_{7}\gamma_{4} - 2c_{6}\gamma_{5} + 6\gamma_{3}\gamma_{4}^2 + c_{8}\gamma_{3}) 
             + t^2(-3c_{7}\gamma_{3} + 3c_{6}\gamma_{4} - 20\gamma_{4}\tilde{\gamma}_{6}) \\
            &+ t^3(-12\gamma_{4}\gamma_{5} + 5c_{6}\gamma_{3})  
             + t^4(-2\gamma_{3}\gamma_{5} - 15\gamma_{4}^2 + 2c_{8}) + t^6(-10c_{6} + 20\tilde{\gamma}_{6}) \\
            &+ 12t^7 \gamma_{5} + 19t^8 \gamma_{4} - 6t^9 \gamma_{3} - 2t^{12}, \\
  \rho_{14} &= c_{7}^2 + 6c_{7}\gamma_{3}\gamma_{4} - 2c_{6}\gamma_{3}\gamma_{5} + 5c_{6}c_{8} 
               - 14c_{8}\tilde{\gamma}_{6} + 6\gamma_{4}\gamma_{10} + 4\gamma_{3}\gamma_{5} \tilde{\gamma}_{6} 
               - 6\gamma_{4}^2 \tilde{\gamma}_{6} \\
            &+ t(-c_{8}\gamma_{5} + 12\gamma_{3}\gamma_{4}\tilde{\gamma}_{6} - 4c_{7}\tilde{\gamma}_{6}) 
             + t^2(-c_{7}\gamma_{5} - 9c_{8}\gamma_{4} - 40\tilde{\gamma}_{6}^2 + 18c_{6}\tilde{\gamma}_{6}) \\
            &+ t^3(-9c_{7}\gamma_{4} + 3c_{6}\gamma_{5} + 10c_{8}\gamma_{3} - 24\gamma_{5} \tilde{\gamma}_{6}) 
             + t^4(-6c_{7}\gamma_{3} - 6\gamma_{10} - 30\gamma_{4} \tilde{\gamma}_{6}) \\
            &+ 12t^5 \gamma_{3} \tilde{\gamma}_{6} - 5t^6 c_{8} + 9t^7 c_{7} + 2t^8 \tilde{\gamma}_{6}, 
\end{align*} 
\begin{align*} 
  \rho_{18} &= \gamma_{9}^2 -2c_{6}c_{7}\gamma_{5} -6c_{7}\gamma_{3}\gamma_{4}^2 + 2c_{7}^2 \gamma_{4} 
             + 2c_{6}\gamma_{3}\gamma_{4}\gamma_{5} - 2c_{6}\gamma_{3}\gamma_{9} 
             - c_{8}\gamma_{4}\tilde{\gamma}_{6}   \\
            &+ 30c_{6} \tilde{\gamma}_{6}^2 + 6c_{7}\gamma_{5} \tilde{\gamma}_{6} - 9c_{6}^2 \tilde{\gamma}_{6} 
             + 6\gamma_{4}^3 \tilde{\gamma}_{6} - 6\gamma_{3}\gamma_{4}\gamma_{5} \tilde{\gamma}_{6} 
             - 20 \tilde{\gamma}_{6}^3 - 6\gamma_{4}^2 \gamma_{10} \\
            &- 9c_{8}\gamma_{10} + 2\gamma_{3}\gamma_{5} \gamma_{10} + c_{6}c_{8}\gamma_{4}
             - 9c_{7}c_{8}\gamma_{3}  \\ 
            &+ t(6c_{7}^2 \gamma_{3} - 24c_{6}c_{7}\gamma_{4} + 7c_{8}\gamma_{4}\gamma_{5} 
             + 6c_{6}\gamma_{5} \tilde{\gamma}_{6} - c_{8}\gamma_{9} + 6\gamma_{3}\gamma_{4}\gamma_{10} 
             - 3c_{8}\gamma_{3} \tilde{\gamma}_{6} \\
            & + 2c_{6}c_{8}\gamma_{3} + 10c_{7}\gamma_{10} + 48c_{7}\gamma_{4} \tilde{\gamma}_{6} 
              - 18\gamma_{3}\gamma_{4}^2 \tilde{\gamma}_{6}) \\
            & + t^2 (25c_{7}\gamma_{4}\gamma_{5} - c_{7}\gamma_{9} + 18c_{6}c_{7}\gamma_{3} 
              - 9c_{6}\gamma_{4} \tilde{\gamma}_{6} - 9c_{8}\gamma_{4}^2 + 60\gamma_{4} \tilde{\gamma}_{6}^2
             - 20 \tilde{\gamma}_{6} \gamma_{10} \\
            & + 9c_{6}\gamma_{10} - 2c_{8}^2 - 31c_{7}\gamma_{3} \tilde{\gamma}_{6} - 4c_{8}\gamma_{3}\gamma_{5}) \\
            & + t^3 (45c_{7}\gamma_{4}^2 - 20c_{7}\gamma_{3}\gamma_{5} - 3c_{6}\gamma_{4}\gamma_{5} 
              + 3c_{6}\gamma_{9} - 15c_{6}\gamma_{3} \tilde{\gamma}_{6} + 11c_{8}\gamma_{3}\gamma_{4}  + 17c_{7}c_{8}  \\
            & + 36\gamma_{4}\gamma_{5} \tilde{\gamma}_{6} - 12\gamma_{5}\gamma_{10}) \\
            & + t^4 (-11c_{7}^2 - 2c_{6}\gamma_{3}\gamma_{5} - 48c_{7}\gamma_{3}\gamma_{4} - 9\gamma_{4}\gamma_{10} 
              - 16c_{8} \tilde{\gamma}_{6} + 5c_{6}c_{8} + 45\gamma_{4}^2 \tilde{\gamma}_{6}  
              + 6\gamma_{3}\gamma_{5} \tilde{\gamma}_{6}) \\ 
            & + t^5 (-51c_{6}c_{7} - 11c_{8}\gamma_{5} + 6\gamma_{3}\gamma_{10} + 120c_{7} \tilde{\gamma}_{6}) \\
            & + t^6 (53c_{7} \gamma_{5} - 60\tilde{\gamma}_{6}^2 + 30c_{6} \tilde{\gamma}_{6} - 10c_{8}\gamma_{4}) \\
            & + t^7 (69c_{7}\gamma_{4} + 3c_{6}\gamma_{5} + c_{8}\gamma_{3} - 36\gamma_{5} \tilde{\gamma}_{6}) \\
            & + t^8 (-16c_{7}\gamma_{3} - 2\gamma_{10} - 57\gamma_{4} \tilde{\gamma}_{6}) \\
            & + 18t^9 \gamma_{3} \tilde{\gamma}_{6} + 2t^{10}c_{8} - 15t^{11} c_{7} + 6t^{12} \tilde{\gamma}_{6}. 
\end{align*}

\begin{rem} 
By using the element $\tilde{\gamma}_{6}$,  it is easy to see that the relations $\rho_{i}' \; (i = 1, 2, 3, 4, 5, 6, 8, 9, 10, 12, 14, 18 )$  in $H^*(E_{7}/T';\Z)$   $($see Theorem $\ref{thm:E_7/T})$  correspond to  the relations $\rho_{i} \; (i = 1, 2, 3, 4, 5, 8, 9, 10, 12, 14, 18), \;  \tilde{\rho}_{6}$   in $H^*(E_{8}/T;\Z)$ 
under the surjective homomorphism $i^{*}: H^{*}(E_{8}/T;\Z) \longrightarrow H^*(E_{7}/T';\Z)$
$($just put $\tilde{\gamma}_{6} =  c_{8} = \gamma_{10} = 0$ in the latter relations$)$.  

\end{rem}


\end{document}